\documentclass[a4paper,10pt,twoside]{amsart}

\usepackage{amssymb}

\def\BC{\mathbb{C}}
\def\BZ{\mathbb{Z}}
\def\BE{\mathbb{E}}

\def\phi{\varphi}
\def\void{\varnothing} %\def\void{\emptyset}

\def\geq{\geqslant}
\def\leq{\leqslant}

\def\bigpr#1{\bigl(#1\bigr)}
\def\Bigbk#1{\Bigl[#1\Bigr]}
\def\condset#1#2{\{\,#1\if#2\empty\else\,;#2\fi\,\}}

\def\Uvoid{U_{\!\void}}
\def\Usubset#1{U_{#1} \setminus \Uvoid}
\def\invzero#1#2#3{(\phi_{#1})^{-1}%
  \bigl(\bigl\{x_{#1}^{#2}\if#3\empty\else=x_{#1}^{#3}\fi=0\bigr\}\bigr)}

\newcommand{\typeAone}[1]{A^{(1)}_{#1}}
\newcommand{\sysAone}[1]{\bigl(\typeAone{#1}\bigr)}
\newcommand{\WtypeAone}[1]{\widetilde{W}\bigl(\typeAone{#1}\bigr)}
\newcommand{\fivealphas}[5]
  {#1\alpha_0 #2\alpha_1 #3\alpha_2 #4\alpha_3 #5\alpha_4}
\DeclareMathOperator{\Res}{Res}

\newcommand{\Peq}[1]{P_\textrm{\uppercase\expandafter{\romannumeral #1}}}

\theoremstyle{plain}     %% caption: bold, body: italic
\newtheorem{thm}{Theorem}[section]

\newtheorem{prop}[thm]{Proposition}
\theoremstyle{definition}%% caption: bold, body: normal

\theoremstyle{remark}    %% caption: italic, body: normal
\newtheorem{rem}[thm]{Remark}

\numberwithin{equation}{section}

\title[An augmentation of the phase space
       of the system of type \protect\(\typeAone{4}\protect\)]
      {An augmentation of the phase space\\%
       of the system of type \protect\(\typeAone{4}\protect\)}
\author{Nobuhiko TAHARA}
\date{}

\begin{document}

\begin{abstract}
We investigate the differential system with affine Weyl group symmetry of
type~\(\typeAone{4}\) and construct a space which parametrizes all
meromorphic solutions of~it.  To demonstrate our method based on
singularity analysis and affine Weyl group symmetry, we first study the
system of type~\(\typeAone{2}\), which is the equivalent of the fourth
Painlev\'e equation, and obtain the space which augments the original phase
space of the system by adding spaces of codimension~1.  For the system of
type~\(\typeAone{4}\), codimension~2 spaces should be added to the phase
space of the system in addition to codimension~1 spaces.
\end{abstract}

\maketitle

\section{Introduction}

The differential system of type~\(\typeAone{l}\) (\(l=2,3,\ldots\,\)),
proposed by~Noumi and~Yamada~\cite{NY}, is a system of autonomous ordinary
differential equations for \((l+1)\) unknown functions \(f_0,\ldots,f_l\)
with complex parameters \(\alpha_0,\ldots ,\alpha_l\) satisfying
\(\alpha_0 + \cdots + \alpha_l = 1\).  The system has the symmetry of the
affine Weyl group of type~\(\typeAone{l}\), where
\(\alpha_0,\ldots,\alpha_l\) are considered as simple roots of the affine
root system of type~\(\typeAone{l}\).  It should be noted that the system
of type~\(\typeAone{l}\) is equivalent to the fourth Painlev\'e equation
when~\(l=2\), and to the fifth when~\(l=3\).

The purpose of this paper is to construct a parameter space of all
meromorphic solutions (including holomorphic solutions) of the system of
type~\(\typeAone{4}\).  We call this space an \emph{augmented} phase space,
and the process to obtain an augmented phase space an \emph{augmentation}
of the original phase space, which is the parameter space of all
holomorphic solutions.  Such spaces have been constructed in the case of
Painlev\'e equations as ``spaces of initial conditions''~\cite{O1}.
They are constructed by means of successive
blowing-up procedures at accessible singular points.  But in our case,
namely for the system of type~\(\typeAone{4}\), the calculations in
blowing-up procedures are very complicated, so we take another approach
based on singularity analysis, that is, the construction of the formal
meromorphic solutions of the system in the form of Laurent expansions at an
arbitrary point containing several arbitrary constants.  We also analyse
the noteworthy connection between the formal meromorphic solutions and the
affine Weyl group symmetry the system has.

In the case of~\(l=2n\) (\(n=1,2,\ldots\,\)), the differential system of
type~\(\typeAone{l}\) is defined by
\[
\sysAone{l}\colon\quad
f_i' = f_i \,\biggl( \sum_{1\leq r\leq n} f_{i+2r-1} -
                     \sum_{1\leq r\leq n} f_{i+2r} \biggr) + \alpha_i,
\]
where \(i = 0,\ldots,2n\) and \('\)~stands for the derivation with respect
to the independent variable~\(t\).  The indices of \(f\) and~\(\alpha\)
should be read as elements of~\(\BZ/(l+1)\BZ\).  The system~\(\sysAone{l}\)
gives a structure of a differential field to the field of rational
functions~\(\BC(\alpha; f)\) of \(\alpha = (\alpha_0,\ldots,\alpha_l)\) and
\(f = (f_0,\ldots,f_l)\), and the system admits an action of the extended
affine Weyl group~\(\WtypeAone{l}\) of type~\(\typeAone{l}\) as B\"acklund
transformation group.  Here B\"acklund transformation means an automorphism
of the differential field~\(\BC(\alpha; f)\) which commutes with the
derivation.  The group~\(\WtypeAone{l}\) is generated by the automorphisms
\(s_0,\ldots ,s_l\) and~\(\pi\) with fundamental relations
\[\begin{array}{ccccc}
  s_i^2 = 1, &
  s_i s_j s_i = s_j s_i s_j & (j = i \pm 1), &
  s_i s_j = s_j s_i & (j \neq i \pm  1); \\
  \pi^{l+1} = 1, &
  \pi s_i = s_{i+1}\pi & & &
\end{array}\]
for \(i,j = 0,1,\ldots,l\).
The actions of \(s_i\)~and~\(\pi\) on \(\alpha_j\)~and~\(f_j\) are given by
\begin{equation}\label{defBT}
  s_i(\alpha_j) = \alpha_j - \alpha_i a_{ij}, \quad
  s_i(f_j) = f_j + \frac{\alpha_i}{f_i} u_{ij}; \quad
  \pi(\alpha_j) = \alpha_{j+1}, \quad
  \pi(f_j) = f_{j+1}
\end{equation}
for \(i,j = 0,1,\ldots,l\), where \(A = (a_{ij})_{0 \leq i, j \leq l}\) is
the generalized Cartan matrix of type~\(\typeAone{l}\):
\begin{equation}\label{defGCM}
  a_{jj} = 2, \quad
  a_{ij} = -1 \ \ (j = i \pm 1), \quad
  a_{ij} = 0  \ \ (j \neq i, i \pm 1),
\end{equation}
and \(U = (u_{ij})_{0 \leq i, j \leq l}\) is the orientation matrix:
\begin{equation}\label{defOM}
  u_{ij} = \pm 1 \ \ (j = i \pm 1), \quad
  u_{ij} = 0 \ \ (j \neq i \pm 1).
\end{equation}

In this paper we investigate the system~\(\sysAone{l}\) in the case of
\(l=2\) and~\(l=4\), namely, the system~\(\sysAone{2}\) and the
system~\(\sysAone{4}\), respectively.
Although the main object of the paper is the system~\(\sysAone{4}\),
we first study the system~\(\sysAone{2}\), in Section~\ref{sec:sys2},
to clarify our method.  Section~\ref{sec:sys2} is divided into
three subsections.  In the first subsection, we obtain all formal
meromorphic solutions of the system~\(\sysAone{2}\), observing where
arbitrary constants appear.  These solutions are classified into three
families of formal meromorphic solutions, each of which corresponds to
B\"acklund transformation~\(s_i\) for some \(i = 0,1,2\).  In the second
subsection, we choose an appropriate coordinate system for each family of
formal meromorphic solutions in order to extract arbitrary constants, or more
precisely, in order that the new coordinates express the formal meromorphic
solutions as formal holomorphic solutions and so that the arbitrary
constants appear in their constant terms.  By the use of such appropriate
coordinate system, the convergence of formal meromorphic solutions is
shown.  In fact, for each family of meromorphic solutions, we obtain two
different coordinate systems.  In the third subsection, we construct a
fiber space~\(\BE\) over the space of parameters~\(\alpha\), so that each
fiber of the space~\(\BE\)
parametrizes all holomorphic solutions and meromorphic solutions.
Holomorphic solutions correspond to the points of a 3-dimensional affine
subspace of the fiber, and three families of meromorphic solutions correspond
to three 2-dimensional affine subspaces.  We also study mappings from~\(\BE\)
to itself associated with B\"acklund transformations.
Our coordinate systems of~\(\BE\) are convenient for studying
these mappings.

The latter sections are devoted to the study of the system~\(\sysAone{4}\).
Sections \ref{sec:sys4-1},~\ref{sec:sys4-2} and~\ref{sec:sys4-3} are
\(\typeAone{4}\)~versions of Subsections \ref{sec:sys2-1},~\ref{sec:sys2-2}
and~\ref{sec:sys2-3}, respectively.  In Section~\ref{sec:sys4-1}, we
study formal meromorphic solutions of the system~\(\sysAone{4}\).  There
are fifteen families of formal meromorphic solutions, which are divided
into three classes corresponding to \(s_i\),~\(s_j s_i\) and~\(s_k s_j s_i\),
or rather the type \((i)\),~\((ij)\) and~\((ijk)\) for some
\(i,j,k=0,\ldots,4\).  In Section~\ref{sec:sys4-2}, we choose a suitable
coordinate system for each family of formal meromorphic solutions.  The
process in the case of type~\((i)\)
is almost the same as the system~\(\sysAone{2}\).
In the case of type~\((ij)\), we first apply the process corresponding to
the type~\((i)\) and then proceed to the next and final process.  In
the case of type~\((ijk)\), we apply the process of the type~\((ij)\), and
then proceed to the final process.  In Section~\ref{sec:sys4-3}, we
construct a fiber space~\(\BE\) over the space of parameters~\(\alpha\) of
the system~\(\sysAone{4}\) in the same way as in
Subsection~\ref{sec:sys2-3}, and study the mappings from~\(\BE\) to itself.
In this case, each fiber~\(\BE(\alpha)\) consists of a 5-dimensional space,
five 4-dimensional spaces and ten 3-dimensional spaces, any two of which do
not intersect.  The 5-dimensional space is a parameter space of the
holomorphic solutions, and each 4-dimensional or 3-dimensional space is a
parameter space of a 4-parameter or a 3-parameter family of meromorphic
solutions, respectively.

The results in Section~\ref{sec:sys4-3} are easily translated to the notion
of the defining manifold (a fiber space over the space of independent
variable~\(t\), whose fibers are the spaces of initial conditions) of the
Hamiltonian system associated with the system~\(\sysAone{4}\).  In the last
section, Section~\ref{sec:sys4-hsys}, we give a list of local coordinate
systems of the defining manifold and the form of the Hamiltonian functions
on the charts.

\section{The system of type \protect\(\typeAone{2}\protect\)}
\label{sec:sys2}

The system of differential equations with affine Weyl group symmetry of
type~\(\typeAone{2}\) is explicitly written as
\begin{equation}\label{sysA12}
\begin{array}{c}
  f_0' = f_0 (f_1 - f_2) + \alpha_0, \\
  f_1' = f_1 (f_2 - f_0) + \alpha_1, \\
  f_2' = f_2 (f_0 - f_1) + \alpha_2,
\end{array}
\end{equation}
where \(' = d/dt\).

\subsection{Formal meromorphic solutions}
\label{sec:sys2-1}

For an arbitrarily fixed~\(t_0\in\BC\), let us consider a formal
meromorphic solution of the system~\eqref{sysA12} of the form
\begin{equation}\label{formalsol2}
f_i = \sum_{n=-r}^\infty c_n^i\, T^n, \quad T := t - t_0 \quad (i = 0,1,2)
\end{equation}
where \(r\)~is a positive integer and
\(c_{-r} = \bigpr{c_{-r}^0, c_{-r}^1, c_{-r}^2} \neq (0,0,0)\).

Substituting this into~\eqref{sysA12} and comparing the coefficients
of~\(T^{n-1}\) on both sides, we have
\begin{equation}\label{syscoeff2}
n c_n^i = \sum_{k=-r}^{n+r-1} c_k^i\, G_{n-k-1}^i + \delta_{0,n-1} \alpha_i
\quad (i = 0,1,2)
\end{equation}
for~\(n \geq -2r+1\) where \(G_n^i = c_n^{i+1} - c_n^{i+2}\) and
\(\delta_{\cdot\,\cdot}\)~is the Kronecker delta.  We set
\(c_n = \bigpr{c_n^0,c_n^1,c_n^2} = (0,0,0)\) for any \(n \leq -r-1\),
by convention.

We first see \(r = 1\) by deriving a contradiction.
For this purpose, we assume \(r > 1\) and
look into the equations~\eqref{syscoeff2} for
\(n = -2r+1, -2r+2, \ldots, -r\):
\begin{align*}
0 &= c_{-r}^i\, G_{-r}^i, \\
0 &= c_{-r}^i\, G_{-r+1}^i + c_{-r+1}^i\, G_{-r}^i, \\
0 &= c_{-r}^i\, G_{-r+2}^i + c_{-r+1}^i\, G_{-r+1}^i + c_{-r+2}^i\, G_{-r}^i,\\
  & \qquad \vdots \\
0 &= c_{-r}^i\, G_{-2}^i + c_{-r+1}^i\, G_{-3}^i
     + \cdots + c_{-3}^i\, G_{-r+1}^i + c_{-2}^i\, G_{-r}^i, \\
 -r c_{-r}^i &= c_{-r}^i\, G_{-1}^i + c_{-r+1}^i\, G_{-2}^i
     + \cdots + c_{-2}^i\, G_{-r+1}^i + c_{-1}^i\, G_{-r}^i.
\end{align*}
The first equation, for~\(i=0,1,2\), is solved as
\[
  c_{-r} = (a,a,a),\ (a,0,0),\ (0,a,0),\ (0,0,a),
\]
where \(a\)~is an arbitrary non-zero constant.
Note that \(G_{-r}^i + c_{-r}^i \neq 0\) (\(i = 0,1,2\)) in every case of the
values of~\(c_{-r}\).
Then these equations enables us to find out either \(c_{-1}^i = 0\) or
\(G_{-1}^i = -r\), for each \(i = 0,1,2\), as follows.
If \(c_{-r}^i = 0\) for some \(i\), then \(G_{-r}^i \neq 0\) by
\(G_{-r}^i + c_{-r}^i \neq 0\) and the above equations yield
\(c_{-r+1}^i = \cdots = c_{-1}^i = 0\).  On the other hand,
if \(c_{-r}^i \neq 0\) for some~\(i\), we obtain \(G_{-r}^i = 0\)
from the first equation, and hence, by the other equations,
\(G_{-r+1}^i = G_{-r+2}^i = \cdots = G_{-2}^i = 0\) and \(G_{-1}^i = -r\).
Therefore \(c_{-r}^i = 0\) for some \(i\) implies \(c_{-1}^i = 0\) and
\(c_{-r}^i \neq 0\) for some~\(i\) implies \(G_{-1}^i = -r\).

Let us derive a contradiction in each case of \(c_{-r}\).
In the case of \(c_{-r} = (0,0,a)\), for~\(a \neq 0\), we obtain
\(c_{-1}^0 = c_{-1}^1 = 0\), \(G_{-1}^2 = -r\) by the above and, since
\(G_{-1}^2 = c_{-1}^0 - c_{-1}^1\) by definition, \(r=0\), which
contradicts the assumption~\(r > 1\).  We have the contradiction in the
case of \(c_{-r} = (0,a,0)\) or \(c_{-r} = (a,0,0)\), similarly.
In the case of \(c_{-r} = (a,a,a)\), for \(a \neq 0\), we have
\(G_{-1}^0 = G_{-1}^1 = G_{-1}^2 = -r\), which contradicts
\(G_{-1}^0 + G_{-1}^1 + G_{-1}^2 =
\bigpr{c_{-1}^1 - c_{-1}^2} + \bigpr{c_{-1}^2 - c_{-1}^0} +
\bigpr{c_{-1}^0 - c_{-1}^1} = 0\).
Thus we have shown that \(r=1\).

Now we determine the coefficients \(c_n = \bigpr{c_n^0,c_n^1,c_n^2}\)
of the expansion~\eqref{formalsol2} for~\(n \geq -1\)
by the equations~\eqref{syscoeff2}.
For~\(n = -1\), we have
\[
  (-1) c_{-1}^i = c_{-1}^i \bigpr{c_{-1}^{i+1} - c_{-1}^{i+2}}
  \quad (i = 0,1,2)
\]
and it follows that
\[
  c_{-1} = (-1,0,1),\ (1,-1,0),\ (0,1,-1).
\]
For~\(n \geq 0\) the equations~\eqref{syscoeff2} can be written as a system
of linear equations
\[
  \bigpr{n - G_{-1}^i} c_n^i - c_{-1}^i c_n^{i+1} + c_{-1}^i c_n^{i+2}
  = \sum_{k=0}^{n-1} c_k^i\, G_{n-k-1}^i + \delta_{0,n-1} \alpha_i
  \quad (i=0,1,2)
\]
with respect to~\(c_n = \bigpr{c_n^0, c_n^1, c_n^2}\), and hence
the coefficients~\(c_n\) are
successively determined as polynomials of
%\(\condset{c_k^i, \alpha_i}{i = 0,1,2,\, k = 0,\ldots,n-1}\),
\(\bigl\{\,c_k^i, \alpha_i\,; i = 0,1,2,\, k = 0,\ldots,n-1\,\bigr\}\)
unless \(\det P_n = 0\),
where \(P_n\) is the coefficient matrix of the linear system:
\[
  P_n =
  \begin{bmatrix}
    n-G_{-1}^0 & -c_{-1}^0 &  c_{-1}^0 \\
      c_{-1}^1 &n-G_{-1}^1 & -c_{-1}^1 \\
     -c_{-1}^2 &  c_{-1}^2 &n-G_{-1}^2
  \end{bmatrix}.
\]
When \(\det P_n = 0\) for some~\(n\) and the linear system has a solution,
the solution contains arbitrary constants, the number of which is equal
to~\(\text{dim}\,\text{ker}\, P_n\).

Let us observe the expansion~\eqref{formalsol2} precisely in the case
of~\(c_{-1} = (-1,0,1)\).  The expansions in the other cases are easily
obtained by the use of cyclic rotations.  We first note that
\[
  P_n =
  \begin{bmatrix}
    n+1& 1 &-1 \\
     0 &n-2& 0 \\
    -1 & 1 &n+1
  \end{bmatrix};\quad
  \det P_n = (n+2) n (n-2).
\]
We can see that \(\text{dim}\,\text{ker}\, P_0 = 1\) and the linear system
for~\(n = 0\) is solved as \(c_0^1 = 0\) and \(c_0^0 = c_0^2\), which is
arbitrary.  The coefficients~\(c_1^i\) for~\(i = 0,1,2\) are
uniquely determined depending on the value of~\(c_0^0 = c_0^2\).
For~\(n = 2\), we can verify that \(\text{dim}\,\text{ker}\, P_2 = 1\) and
the solution \(c_2 = \bigpr{c_2^0, c_2^1, c_2^2}\) of the linear system is
determined so that \(c_2^1\)~is an arbitrary constant while the other
\(c_2^0\) and~\(c_2^2\) are unique depending on the value of~\(c_2^1\).
The coefficients~\(c_n\) for~\(n \geq 3\) are uniquely determined by
\(c_0^0\) and~\(c_2^1\).  Thus we have obtained the expansion
of the formal meromorphic solution as
\begin{equation}\label{series2type1}
\begin{split}
f_0 &= -\frac{1}{T} + c_0^0
       + \frac{(2\alpha_0 + 3\alpha_1 + \alpha_2) - \bigpr{c_0^0}^2}{3}\, T
       + O(T^2),\\
f_1 &= \phantom{-\frac{1}{T} + c_0^0}-\alpha_1 T + c_2^1\, T^2 + O(T^3), \\
f_2 &= \phantom{-}\frac{1}{T} + c_0^0
       + \frac{(\alpha_0 + 3\alpha_1 + 2\alpha_2) + \bigpr{c_0^0}^2}{3}\, T
       + O(T^2),
\end{split}
\end{equation}
where \(c_0^0\) and~\(c_2^1\) are arbitrary constants which are free from
the system~\eqref{sysA12}.  Hence this formal solution depends on two
arbitrary constants.  Note that we do not
consider the position of the pole, \(t_0\), as an arbitrary constant, since
the system~\eqref{sysA12} is autonomous.

Let us set
\[
  \Res f = (\Res_{t=t_0} f_0,\Res_{t=t_0} f_1,\Res_{t=t_0} f_2)
\]
for a set of formal meromorphic functions~\(f = (f_0,f_1,f_2)\)
at~\(t=t_0\).  Then the above solution~\eqref{series2type1} can be
indicated by~\(\Res f = (-1,0,1)\).
All formal solutions of the system~\eqref{sysA12} with \(\Res f = (-1,0,1)\)
defines a 2-parameter family of formal solutions.
The other 2-parameter families are also defined by the formal solutions
with~\(\Res f = (1,-1,0)\) or with~\(\Res f = (0,1,-1)\).

In the rest of this subsection, we mention the relations between these
three 2-parameter families of formal meromorphic solutions and the
B\"acklund transformations~\(s_i\) for~\(i=0,1,2\).
For example, let \(f\)~be the solution with~\(\Res f = (-1,0,1)\) and let
\(g_i = s_1(f_i)\) and \(\beta_i = s_1(\alpha_i)\) for~\(i = 0,1,2\),
namely, from \eqref{defBT},~\eqref{defGCM} and~\eqref{defOM},
\[
    g_0 = f_0 - \frac{\alpha_1}{f_1},
\ \ g_1 = f_1,
\ \ g_2 = f_2 + \frac{\alpha_1}{f_1}; \quad
    \beta_0 = \alpha_0 + \alpha_1,
\ \ \beta_1 = \alpha_1,
\ \ \beta_2 = \alpha_2 + \alpha_1.
\]
Then \(g = (g_0,g_1,g_2)\) satisfies the system
\[
  g_i' = g_i (g_{i+1} - g_{i+2}) + \beta_i, \quad (i=0,1,2)
\]
and it has the series expansion given by
\[
  g_0 = \biggl( c_0^0 + \frac{c_2^1}{\beta_1} \biggr) + O(T), \quad
  g_1 = 0 - \beta_1 T + c_2^1\, T^2 + O(T^3), \quad
  g_2 = \biggl( c_0^0 - \frac{c_2^1}{\beta_1} \biggr) + O(T)
\]
without terms of negative powers of~\(T\).  This fact means that,
when \(\alpha_1 (= \beta_1) \neq 0\), the formal solution with
\(\Res f = (-1,0,1)\) is related, via B\"acklund transformation~\(s_1\),
to the formal holomorphic solution of the system~\eqref{sysA12} with
parameters \(s_1(\alpha) = \beta = (\beta_0,\beta_1,\beta_2)\) under the
initial condition
\[
  g_0(t_0) = c_0^0 + \frac{c_2^1}{\beta_1},\quad
  g_1(t_0) = 0, \quad
  g_2(t_0) = c_0^0 - \frac{c_2^1}{\beta_1}.
\]
In this sense, we say that the solution with \(\Res f = (-1,0,1)\)
corresponds to the B\"acklund transformation~\(s_1\) and that the formal
solution is of type~\((1)\).  Then the solutions with \(\Res f = (1,-1,0)\)
and \((0,1,-1)\) are of type~\((2)\) and~\((0)\), respectively.
Generalizing the terminology,
we say that the 3-parameter family of holomorphic solutions is of type
\((\void)\) (see Table~\ref{tab:res2}).

\begin{table}
\caption{Classification of the families of solutions
         of the system \protect\(\sysAone{2}\protect\)}
\label{tab:res2}
\[\begin{array}{c|c|c|c}
\hline
\text{type}&\Res f &\text{corresponding BT}
  &\text{\# of arbitrary constants} \\ \hline

(\void)& (0,0,0)  & \text{id} & 3 \\ \hline

(1)    & (-1,0,1) &       s_1 &   \\
(2)    & (1,-1,0) &       s_2 & 2 \\
(0)    & (0,1,-1) &       s_0 &   \\
\hline
\end{array}\]
\end{table}

\subsection{Coordinates for formal meromorphic solutions}
\label{sec:sys2-2}

In this subsection, we prove the convergence of the formal meromorphic
solutions~\eqref{series2type1} by choosing suitable coordinate systems so
that the arbitrary constants in the formal solutions are interpreted as
initial conditions of the holomorphic differential systems in the new
coordinate systems.

Let \(f = (f_0, f_1, f_2)\) be the solution~\eqref{series2type1}
with~\(\Res f = (-1,0,1)\).  We first notice that
\(f_0 + f_1 + f_2 = 2 c_0^0 + T\) and that \(f_0\)~has a simple pole while
\(f_1\)~has a simple zero at~\(t = t_0\).  Then \(1/f_0\)~has a simple
zero, and moreover, the expansion of \(f_0 f_1\)~is of the form
\[ f_0 f_1 = \alpha_1 - \bigpr{c_2^1 + \alpha_1 c_0^0}\,T + O(T^2). \]
Noting the position of the arbitrary constant~\(c_2^1\),
we multiply \(\alpha_1 - f_0 f_1\) by~\(f_0\).  Then we have
\[ f_0 (\alpha_1 - f_0 f_1) = \bigpr{c_2^1 + \alpha_1 c_0^0} + O(T), \]
where \(c_2^1\)~is contained in the initial value at~\(t=t_0\).

Based on such an observation, we introduce the transformation
\begin{equation}\label{patch2type1m}
  u_0 = 1/f_0, \quad
  u_1 = f_0 (\alpha_1 - f_0 f_1), \quad
  u_0 + u_1 + u_2 = f_0 + f_1 + f_2.
\end{equation}
We can verify that the transformation~\eqref{patch2type1m} changes
the system~\eqref{sysA12} into the holomorphic system
\begin{equation}\label{sysA12_1}\begin{aligned}
  u_0' &=
    2 u_0^3 u_1 - (\alpha_0 + 2\alpha_1 - 1) u_0^2
    + (u_1 + u_2) u_0 - 1, \\
  u_1' &=
    - 3 u_0^2 u_1^2 + (2\alpha_0 + 4\alpha_1 - 1) u_0 u_1
    - (u_1 + u_2) u_1 - (\alpha_0 + \alpha_1) \alpha_1, \\
  u_0' + u_1' + u_2' &= \alpha_0 + \alpha_1 + \alpha_2.
\end{aligned}\end{equation}
Now take a solution of the new system~\eqref{sysA12_1} holomorphic
at~\(t = t_0\) with \(u_0(t_0) = 0\), \(u_1(t_0) = h_1\), \(u_2(t_0) = h_2\)
by Cauchy's existence and uniqueness theorem.  Then, since
\[
  u_0'(t_0) = -1,\quad u_0''(t_0) = -(h_1 + h_2),
\]
we have the Taylor expansion of the holomorphic solution as
\begin{align*}
  u_0 &= -T - \tfrac{1}{2} (h_1 + h_2)\,T^2 + O(T^3), \\
  u_1 &= h_1 + O(T), \\
  u_2 &= h_2 + O(T).
\end{align*}
Transform this solution~\(u(t) = (u_0(t),u_1(t),u_2(t))\)
by~\eqref{patch2type1m}, then we have a meromorphic solution
of~\eqref{sysA12} expanded as
\begin{align*}
  f_0 &= 1 / u_0
       = -\tfrac{1}{T} + \tfrac{1}{2} (h_1 + h_2) + O(T), \\
  f_1 &= u_0 (\alpha_1 - u_0 u_1)
       = -\alpha_1 T - \bigl(h_1 + \tfrac{1}{2}\alpha_1(h_1+h_2)\bigr)\,T^2
         + O(T^3), \\
  f_2 &= u_0 + u_1 + u_2 - 1/u_0 - u_0 (\alpha_1 - u_0 u_1)
       = \tfrac{1}{T} + \tfrac{1}{2} (h_1 + h_2) + O(T),
\end{align*}
which coincides formally with the formal solutions~\eqref{series2type1}
under the one-to-one correspondence
\[
  c_0^0 = \tfrac{1}{2} (h_1 + h_2), \quad
  c_2^1 = h_1 + \tfrac{1}{2} \alpha_1 (h_1 + h_2).
\]
Thus we have proved the convergence of the formal meromorphic
solution~\eqref{series2type1} with~\(\Res f = (-1,0,1)\), or of type~\((1)\).

From the above result, it follows that the domain of definition of the
system~\eqref{sysA12} can be extended from the affine space
\(\condset{(f_0,f_1,f_2)\in\BC^3}{}\) to the space obtained by
identification of \(\condset{(f_0,f_1,f_2)\in\BC^3}{}\) and
\(\condset{(u_0,u_1,u_2)\in\BC^3}{}\) via~\eqref{patch2type1m}.
The identified space is considered to be a disjoint union of the original
phase space \(\condset{(f_0, f_1, f_2)\in\BC^3}{}\) and the 2-dimensional
affine space \(\condset{(u_0,u_1,u_2)\in\BC^3}{u_0 = 0}\).  The added
2-dimensional affine space is considered to be a parameter space of the
2-parameter family of solutions with~\(\Res f = (-1,0,1)\).  We call such
an extension of the domain of definition an \emph{augmentation} of the
phase space.

We note that the same argument can be done by another transformation
\begin{equation}\label{patch2type1p}
  v_2 = 1/f_2, \quad
  v_1 = f_2 (-\alpha_1 - f_2 f_1), \quad
  v_0 + v_1 + v_2 = f_0 + f_1 + f_2.
\end{equation}
The affine space \(\condset{(v_0,v_1,v_2)\in\BC^3}{v_2 = 0}\) is also the
parameter space of the same 2-parameter family of solutions
with~\(\Res f = (-1,0,1)\) and isomorphic to \(\{u_0 = 0\}\).

In the end of this subsection, we observe how \(s_i\) (\(i=0,1,2\)) act on
the variables \(u_0, u_1, u_2\) and \(v_0,v_1,v_2\).
The actions of~\(s_0\) on~\(u_0, u_1, u_2\) are calculated as
\begin{align*}
  s_0(u_0) &= s_0(1/f_0) = 1 / s_0(f_0) = 1/f_0 = u_0, \\
  s_0(u_1) &= s_0(f_0 (\alpha_1 - f_0 f_1)) \\
           &= s_0(f_0) (s_0(\alpha_1) - s_0(f_0) s_0(f_1) ) \\
           &= f_0 \Bigl(\alpha_0 + \alpha_1 - f_0
                        \bigl(f_1 + \tfrac{\alpha_0}{f_0}\bigr)\Bigr) \\
           &= f_0 (\alpha_1 - f_0 f_1) = u_1, \\
  s_0(u_2) &= (f_0 + f_1 + f_2) - (u_0 + u_1) = u_2.
\end{align*}
Similarly, we obtain the action of \(s_1\) as
\[
  s_1(u_0) = u_0 - \frac{\alpha_1}{u_1}, \quad
  s_1(u_1) = u_1, \quad
  s_1(u_2) = u_2 + \frac{\alpha_1}{u_1}.
\]
However we see that the forms of~\(s_2(u_i)\) for~\(i = 0,1,2\) are not so
simple as above.  On the other hand, by~\eqref{patch2type1p}, we see
\begin{align*}
  s_1(v_0) = v_0 - \frac{\alpha_1}{v_1}, \quad
  s_1(v_1) &= v_1, \quad
  s_1(v_2) = v_2 + \frac{\alpha_1}{v_1}, \\
  s_2(v_0) = v_0, \quad s_2(v_1) &= v_1, \quad s_2(v_2) = v_2,
\end{align*}
but the forms of~\(s_0(v_i)\) for~\(i = 0,1,2\) are complicated.
Hence the coordinate system~\(u\) is convenient to observe the action
of~\(s_0\) while \(v\) is convenient to see the action of \(s_2\).
This is the reason why we take both systems \(u\) and~\(v\).
In the next subsection, \(u\)~and~\(v\) will be distinguished by the labels
with \(-\)~and~\(+\), respectively.

We can obtain the similar results in the cases of
\(\Res f = (1,-1,0)\) and of \(\Res f = (0,1,-1)\), obviously.

\subsection{Augmentation of the phase space and B\"acklund transformations}
\label{sec:sys2-3}

We now define a fiber space~\(\BE\) over the parameter space
\[
  V = \condset{\alpha = (\alpha_0,\alpha_1,\alpha_2) \in \BC^3}
              {\alpha_0 + \alpha_1 + \alpha_2 = 1},
\]
each fiber~\(\BE(\alpha)\) of which is the augmented phase space
parametrizing all meromorphic solutions of the system~\eqref{sysA12}.

Let \(I = \bigl\{\void, 0_+, 1_+, 2_+, 0_-, 1_-, 2_- \bigr\}\) be a label set
and \(W_*\) for each \(* \in I\) be seven copies of~\(V \times \BC^3\)
with coordinates
\((\alpha, x_*)=\bigpr{\alpha_0,\alpha_1,\alpha_2; x_*^0,x_*^1,x_*^2}\in W_*\).
Then we define the space~\(\BE\) by gluing \(W_*\) via the following
identification equations
\[
  x_{i_+}^{i+1} = 1/x_\void^{i+1}, \quad
  x_{i_+}^i = x_\void^{i+1} \bigpr{-\alpha_i - x_\void^{i+1} x_\void^i}, \quad
  x_{i_+}^{i-1} + x_{i_+}^i + x_{i_+}^{i+1}
    = x_\void^{i-1} + x_\void^i + x_\void^{i+1}
\]
for \((\alpha, x_\void) \in W_\void\) and \((\alpha, x_{i_+}) \in W_{i_+}\),
and
\[
  x_{i_-}^{i-1} = 1/x_\void^{i-1}, \quad
  x_{i_-}^i = x_\void^{i-1} \bigpr{+\alpha_i - x_\void^{i-1} x_\void^i}, \quad
  x_{i_-}^{i-1} + x_{i_-}^i + x_{i_-}^{i+1}
    = x_\void^{i-1} + x_\void^i + x_\void^{i+1}
\]
for \((\alpha, x_\void)\in W_\void\) and \((\alpha, x_{i_-})\in W_{i_-}\),
namely,
\[
  \BE=\biggl(\bigsqcup_{*\in I}W_*\biggr)\,\bigg/\sim,
\]
where \(\sim\)~is the equivalence relation generated by the above
equations.  Here \(x_\void\)~shall be considered as the original coordinate
system~\(f\) in~\eqref{sysA12}.  Then we see that \(x_{1_-} = u\) and
\(x_{1_+} = v\), where \(u\)~and~\(v\) are the coordinate systems defined
by~\eqref{patch2type1m} and~\eqref{patch2type1p}, respectively.
We denote by~\(\pi_V\) the natural projection~\(\pi_V \colon \BE \to V\)
and let \(\BE(\alpha) := \pi_V^{-1}(\alpha)\)~be the fiber of~\(\BE\)
over~\(\alpha \in V\).

For \((\alpha, x_*)\in W_*\) (\(* \in I\)),
let \([(\alpha, x_*)]\) be its equivalence class.
Define the subset \(U_*\) of \(\BE\) by
\[U_* := \condset{[(\alpha, x_*)]}{(\alpha, x_*) \in W_*}.\]
Then there are coordinate mappings \(\phi_*\) from~\(U_*\) to~\(W_*\), namely,
\[
  \phi_* \colon p \in U_* \mapsto
  \bigpr{\alpha_0(p), \alpha_1(p), \alpha_2(p); x_*^0(p), x_*^1(p), x_*^2(p)}
  \in W_* = V\times\BC^3
\]
for each \(*\in I\).  Note that \(\alpha_i\)~and~\(x_*^i\) are here
considered to be coordinate functions.  Let
\(\BE_\void\) and \(\BE_i\) (\(i = 0,1,2\)) be the subsets of \(\BE\)
defined by
\begin{align*}
  \BE_\void := \Uvoid,\quad
  \BE_i := \Usubset{i_+} &= \invzero{i_+}{i+1}{} \\
         = \Usubset{i_-} &= \invzero{i_-}{i-1}{}.
\end{align*}
Then
\(\BE_\void \cong V \times \BC^3\), \(\BE_i \cong V \times \BC^2\)
and the space~\(\BE\) is decomposed as
\[
  \BE = \BE_\void \sqcup \BE_0 \sqcup \BE_1 \sqcup \BE_2.
\]
Each \(\BE_i\) (\(i=0,1,2\)) is the parameter space of the
corresponding 2-parameter family of meromorphic solutions of type~\((i)\)
and \(\BE_\void\)~is that of the 3-parameter family of holomorphic
solutions of type~\((\void)\) of the system~\eqref{sysA12}.

The system of differential equations~\eqref{sysA12} defines a vector field
\[
  X_p = \sum_{i=0}^2
        \Bigl(x_\void^i(p) \bigl(x_\void^{i+1}(p) - x_\void^{i+2}(p)\bigr)
              + \alpha_i(p)\Bigr)
        \left(\frac\partial{\partial x_\void^i}\right)_p,
        \quad p\in\BE(\alpha)\cap \Uvoid
\]
for each~\(\alpha\in V\).  As we have already shown, the vector field can
be holomorphically extended to~\(\BE(\alpha)\).

We next observe how B\"acklund transformations act on the space~\(\BE\).
For~\(w\in\WtypeAone{2}\), we define \(\sigma_w\) from~\(\BE\) to itself by
\[
  (\alpha_j\circ\sigma_w)(p) = (w(\alpha_j))(p), \quad
  \bigpr{x_*^j\circ\sigma_w}(p) = \bigpr{w\bigpr{x_*^j}}(p) \quad
    (p\in\BE,\, j = 0,1,2)
\]
for any \(*\in I\) as far as the right-hand sides are defined for~\(p\).
Here \(w\bigpr{x_{*}^j}\) (\(j = 0,1,2\)) are rational functions of
\(\alpha = (\alpha_0,\alpha_1,\alpha_2)\) and \(x_* = (x_*^0, x_*^1, x_*^2)\)
determined by \(w \colon \BC(\alpha, f) \to \BC(\alpha, f)\) and the
isomorphism from \(\BC(\alpha, f)\) to \(\BC(\alpha, x_*)\).

We can verify that \(\sigma_w\) is birational and
\(\sigma_{w'}\circ\sigma_w = \sigma_{ww'}\) for~\(w, w'\in\WtypeAone{2}\).
Note that \(\WtypeAone{2}\)~is generated by \(s_0, s_1, s_2\) and~\(\pi\),
and that \(\sigma_\pi\)~is obviously extended to a biholomorphic mapping
from~\(\BE(\alpha)\) to~\(\BE(\pi(\alpha))\),
which maps \(\BE_*\) (\(* = \void,0,1,2\)) as
\[
  \sigma_\pi(\BE_\void) = \BE_\void, \quad
  \sigma_\pi(\BE_0) = \BE_1, \quad
  \sigma_\pi(\BE_1) = \BE_2, \quad
  \sigma_\pi(\BE_2) = \BE_0.
\]
Therefore we now consider \(\sigma_i:=\sigma_{s_i}\) for \(i = 0,1,2\) in
more detail.  Let \(i = 0,1,2\)~be fixed.  We can verify that
\(\sigma_i\)~can be extended to a biholomorphic mapping
from~\(\BE(\alpha)\) to~\(\BE(s_i(\alpha))\)
for any~\(\alpha\in V\) and its images of~\(\BE_*\)
(\(* = \void,0,1,2\)) are given by
\[
  \begin{array}{ccc}
    \sigma_i(\BE_\void \setminus D_i) = \BE_\void \setminus D_i, &
    \sigma_i(\BE_\void' \cap D_i) = \BE_i', &
    \sigma_i(\BE_{i \pm 1}) = \BE_{i \pm 1} \\
    &
    \sigma_i(\BE_i') = \BE_\void' \cap D_i, &
    \\
  \end{array}
\]
where
\[
  \BE_*' := \BE_* \setminus \pi_V^{-1}(\{\alpha_i = 0\}), \quad
   D_i   := \invzero{\void}{i}{}.
\]
These assertions are verified by direct calculation using
\[
  s_i(x_*^j) =
    \begin{cases}
      x_*^j + \dfrac{\alpha_i}{x_*^i} u_{ij} & (* = \void, i_+, i_-), \\
      x_*^j                                  & \bigpr{* = (i-1)_+, (i+1)_-},
    \end{cases}
\]
where \(u_{ij}\)~are given by~\eqref{defOM}.

We remark that such relations, say
\(\sigma_1(\BE_1') = \BE_\void' \cap D_1\), explain the correspondence of
the meromorphic solutions of type~\((1)\) and the holomorphic solutions of
type~\((\void)\).  \(\sigma_1\)~maps the integral curve through a point
of~\(\BE_1'\) to an integral curve through a point of~\(\BE_\void' \cap D_1\).

\section{Formal meromorphic solutions of the system
         \protect\(\sysAone{4}\protect\)}
\label{sec:sys4-1}

In this section, we obtain all formal meromorphic solutions of the system
of differential equations
with affine Weyl group symmetry of type~\(\typeAone{4}\):
\begin{equation}\label{sysA14}
  f_i' = f_i (f_{i+1} - f_{i+2} + f_{i+3} - f_{i+4}) + \alpha_i
  \quad (i = 0,\ldots,4).
\end{equation}
We show that (i)~the order of the pole of every formal solution is one,
(ii)~there are fifteen families of such formal solutions, and
(iii)~five of them contain four arbitrary constants and the other ten
contain three arbitrary constants.

We proceed in the same manner as Subsection~\ref{sec:sys2-1} of the case
of type~\(\typeAone{2}\).  Let \(f = (f_0,\dots,f_4)\)~be a formal solution
of the form
\begin{equation}\label{formalsol}
  f_i = \sum_{n=-r}^\infty c_n^i\, T^n, \quad T := t - t_0 \quad (i=0,\ldots,4)
\end{equation}
where \(r\)~is a positive integer.
Direct substitution of these series into~\eqref{sysA14} gives
\begin{equation}\label{syscoeff}
  n c_n^i = \sum_{k=-r}^{n+r-1} c_k^i\, G_{n-k-1}^i + \delta_{0,n-1} \alpha_i
  \quad (i=0,\ldots,4)
\end{equation}
where
\[
  G_n^i = \sum_{k=1}^4 (-1)^{k-1} c_n^{i+k}
        = c_n^{i+1} - c_n^{i+2} + c_n^{i+3} - c_n^{i+4}.
\]

\begin{prop}
If \(f_i\) has a pole at \(t = t_0\) then the order of the pole is one,
i.e., \(r = 1\).
\end{prop}

\begin{proof}
Suppose that \(r > 1\).  Consider a series of equations~\eqref{syscoeff}
for~\(n = -2r+1, -2r+2, \ldots,-r\).  The first equation
\(0 = c_{-r}^i\, G_{-r}^i\) (\(i = 0,\ldots,4\)) has the following
solutions up to cyclic rotation
\[
  c_{-r} = \bigpr{c_{-r}^0,\ldots,c_{-r}^4}
         = (a,0,0,0,0),\ (a,a,a,0,0),\ (a,a,0,-a,0),\ (a,a,a,a,a)
\]
where \(a \neq 0\).  It is verified that \(G_{-r}^i + c_{-r}^i \neq 0\) for
each \(i = 0,\ldots,4\) and each solution~\(c_{-r}\).
Therefore, as in Subsection~\ref{sec:sys2-1},
the remaining equations yield that \(c_{-r}^i = 0\) implies
\(c_{-1}^i = 0\) and \(c_{-r}^i \neq 0\) implies \(G_{-1}^i = -r\) for each
\(i = 0,\ldots,4\).  Then, in each case of~\(c_{-r}\),
we derive the contradiction that \(r = 0\) as follows:
% adjust position of label
\begin{list}{(\roman{enumi})\hfill}{\usecounter{enumi}%
  \labelwidth\leftmargini \addtolength\labelwidth{-\labelsep}}
\item In the case of \(c_{-r} = (a,0,0,0,0)\),
  we deduce that
  \(c_{-1}^1 = \cdots = c_{-1}^4 = 0\) and \(G_{-1}^0 = -r\).
  Then it follows that \(r = 0\) from
  \(G_{-1}^0 = c_{-1}^1 - c_{-1}^2 + c_{-1}^3 - c_{-1}^4\).
\item In the case of \(c_{-r} = (a,a,a,0,0)\), we have
  \(G_{-1}^0 = G_{-1}^1 = G_{-1}^2 = -r\), \(c_{-1}^3 = c_{-1}^4 = 0\) and
  \(G_{-1}^0 + G_{-1}^1 + G_{-1}^2 = c_{-r}^3 + c_{-r}^4\),
  and hence \(r = 0\).
\item In the case of \(c_{-r} = (a,a,0,-a,0)\), we have
  \(G_{-1}^0 = G_{-1}^1 = G_{-1}^3 = -r\), \(c_{-1}^2 = c_{-1}^4 = 0\) and
  \(G_{-1}^0 + G_{-1}^1 - G_{-1}^3 = c_{-r}^2 - c_{-r}^4\),
  and hence \(r = 0\).
\item In the case of \(c_{-r} = (a,a,a,a,a)\), we have
  \(G_{-1}^0 = \cdots = G_{-1}^4 = -r\) and
  \(G_{-1}^0 + \cdots + G_{-1}^4
    = \sum_{i=0}^4\sum_{k=1}^4 (-1)^{k-1} c_n^{i+k} = 0\),
  and hence \(r = 0\). \qed
\end{list}
\renewcommand{\qed}{}
\end{proof}

Let us now determine the coefficients~\(c_n\) (\(n \geq -1\)) of the
expansion~\eqref{formalsol} from~\eqref{syscoeff}.
For~\(n = -1\), the equations~\eqref{syscoeff} are written by
\[
  (-1) c_{-1}^i =
  c_{-1}^i \bigpr{c_{-1}^{i+1} - c_{-1}^{i+2} + c_{-1}^{i+3} - c_{-1}^{i+4}}
  \quad (i = 0,\ldots,4)
\]
which has fifteen solutions \(c_{-1} = \bigpr{c_{-1}^0,\ldots,c_{-1}^4}\),
each of which equals to one of
\[
(-1,0,1,0,0),\ (-1,0,0,0,1),\ (-1,-3,0,3,1),
\]
by suitable cyclic rotations.

For~\(n \geq 0\), the equations~\eqref{syscoeff} with respect to~\(c_n\)
are written by the following linear system
\[
  (n - G_{-1}^i) c_n^i
    - c_{-1}^i \bigpr{c_n^{i+1} - c_n^{i+2} + c_n^{i+3} - c_n^{i+4}}
  = \sum_{k=0}^{n-1} c_k^i\, G_{n-k-1}^i + \delta_{0,n-1} \alpha_i.
  \quad (i = 0,\ldots,4)
\]
Let
\[
  P_n := \begin{bmatrix}
           n-G_{-1}^0 & -c_{-1}^0 &  c_{-1}^0 & -c_{-1}^0 &  c_{-1}^0 \\
             c_{-1}^1 &n-G_{-1}^1 & -c_{-1}^1 &  c_{-1}^1 & -c_{-1}^1 \\
            -c_{-1}^2 &  c_{-1}^2 &n-G_{-1}^2 & -c_{-1}^2 &  c_{-1}^2 \\
             c_{-1}^3 & -c_{-1}^3 &  c_{-1}^3 &n-G_{-1}^3 & -c_{-1}^3 \\
            -c_{-1}^4 &  c_{-1}^4 & -c_{-1}^4 &  c_{-1}^4 &n-G_{-1}^4
         \end{bmatrix}.
\]
Then \(c_n\)~is uniquely determined
by~\(c_1,\ldots,c_{n-1}\) unless \(\det P_n = 0\).
Since \(P_n\) only depends on~\(n\) and~\(c_{-1}\),
it suffices to consider the above three typical cases
of the values of~\(c_{-1}\).

% enumerated paragraphs
\begin{list}{(\roman{enumi})}{\usecounter{enumi}%
  \leftmargin0pt \itemindent0pt \labelwidth-\parindent%
  \addtolength{\labelwidth}{-\labelsep}}
\item The case of \(c_{-1} = (-1,0,1,0,0)\):
We have
\[
  \det P_n = (n+2) n^3 (n-2)
 \]
and we have the relations \(c_0^1 = 0\), \(c_0^4 - c_0^0 + c_0^2 - c_0^3 = 0\)
among \(c_0^0,\ldots,c_0^4\).  In this case the system has a formal
meromorphic solution of the form
\begin{equation}\begin{aligned}\label{seriestype1}
  f_0 &=
    - \frac{1}{T} + c_0^0
    + \frac{1}{3}
      \Bigbk{\fivealphas{2}{+3}{+}{-}{+}
             - \bigpr{c_0^0}^2 + 2\bigpr{c_0^2}^2 - 2c^2}\,T
    + O(T^2),\\
  f_1 &= \phantom{-\frac{1}{T} + c_0^0}
    -\alpha_1 T + c_2^1\,T^2 + O(T^3),\\
  f_2 &= \phantom{-}
    \frac{1}{T} + c_0^2
    + \frac{1}{3}
      \Bigbk{\fivealphas{}{+3}{+2}{+}{-}
             + \bigpr{c_0^2}^2 - 2\bigpr{c_0^0}^2 + 2c^2}\,T
    + O(T^2),\\
  f_3 &= \phantom{-\frac{1}{T}+{}}
      \bigpr{c_0^2 + c}
      + \Bigbk{\alpha_3 - \bigpr{c_0^2}^2 + c^2}\,T + O(T^2), \\
  f_4 &= \phantom{-\frac{1}{T}+{}}
      \bigpr{c_0^0 + c}
      + \Bigbk{\alpha_4 + \bigpr{c_0^0}^2 - c^2}\,T + O(T^2),
\end{aligned}\end{equation}
where \(c := c_0^4 - c_0^0 = c_0^3 - c_0^2\).
This formal solution depends on four arbitrary constants
\(c_0^0, c_0^2, c_0^3, c_2^1\)
and therefore defines a 4-parameter family of formal meromorphic solutions.

\item The case of \(c_{-1} = (-1,0,0,0,1)\):
We have
\[
  \det P_n = (n+2)^2 n (n-2)^2
\]
and obtain a formal solution containing three arbitrary constants
\(c_0^0, c_2^1, c_2^3\) written by
\begin{equation}\begin{aligned}\label{seriestype13}
  f_0 &=
    -\frac{1}{T} + c_0^0
    + \frac{1}{3}\Bigbk{\fivealphas{2}{+3}{+}{+3}{+} - \bigpr{c_0^0}^2}\,T
    + O(T^2),\\
  f_1 &= \phantom{-\frac{1}{T}+c_0^0}
    -\alpha_1 T + c_2^1\,T^2 + O(T^3),\\
  f_2 &= \phantom{-\frac{1}{T}+c_0^0}\ \,
    \frac{1}{3}\alpha_2 T \phantom{{}+ c_2^1 T^2}{}+ O(T^3),\\
  f_3 &= \phantom{-\frac{1}{T}+c_0^0}
    -\alpha_3 T + c_2^3\,T^2 + O(T^3),\\
  f_4 &=  \phantom{-}
    \frac{1}{T} + c_0^0
    + \frac{1}{3}\Bigbk{\fivealphas{}{+3}{+}{+3}{+2} + \bigpr{c_0^0}^2}\,T
    + O(T^2).
\end{aligned}\end{equation}
This is a 3-parameter family of formal meromorphic solutions.

\item The case of \(c_{-1} = (-1,-3,0,3,1)\):
We have
\[
  \det P_n = (n+4) (n+2) n (n-2) (n-4)
\]
and get a formal solution including three arbitrary constants
\(c_0^0, c_2^3, c_4^2\) written by
\begin{equation}\begin{aligned}\label{seriestype132}
f_0 &=
  -\frac{1}{T} + c_0^0
  + \frac{1}{3}\Bigbk{\fivealphas{2}{+3}{+5}{+3}{+} - \bigpr{c_0^0}^2}\,T
  - c_2^3\,T^2 + O(T^3),\\
f_1 &=
  -\frac{3}{T} \phantom{{}+c_0^0}
  + \frac{1}{5}\Bigbk{\fivealphas{}{-2}{-5}{-3}{-} - 2\bigpr{c_0^0}^2}\,T \\
    & \phantom{=-\frac{1}{T} + c_0^0+ \frac{1}{3}}\quad
  + \Bigbk{c_2^3 - \frac{1}{2}(\fivealphas{}{+3}{+5}{+3}{+})c_0^0}\,T^2
  + O(T^3),\\
f_2 &= \phantom{-\frac{1}{T}+c_0^0}
  -\frac{1}{3}\alpha_2 T
  + \frac{1}{45}\alpha_2
    \Bigbk{-\alpha_0-3\alpha_1+3\alpha_3+\alpha_4 + 2\bigpr{c_0^0}^2}\,T^3
  + c_4^2\,T^4 + O(T^5),\\
f_3 &= \phantom{-}
  \frac{3}{T} \phantom{{}+c_0^0}
  - \frac{1}{5}\Bigbk{\fivealphas{}{+3}{+5}{+2}{-} - 2\bigpr{c_0^0}^2}\,T
  + c_2^3\,T^2 + O(T^3), \\
f_4 &= \phantom{-}
  \frac{1}{T} + c_0^0
  + \frac{1}{3}\Bigbk{\fivealphas{}{+3}{+5}{+3}{+2} + \bigpr{c_0^0}^2}\,T \\
    & \phantom{=-\frac{1}{T} + c_0^0+ \frac{1}{3}}\quad
  - \Bigbk{c_2^3 - \frac{1}{2}(\fivealphas{}{+3}{+5}{+3}{+})c_0^0}\,T^2
  + O(T^3).
\end{aligned}\end{equation}
This is another 3-parameter family of solutions.
\end{list}

As the system of type~\(\typeAone{2}\), we denote the coefficients
\(c_{-1} = \bigpr{c_{-1}^0,\ldots,c_{-1}^4}\) by
\(\Res f = (\Res f_0,\ldots,\Res f_4)\).
The above results are stated as:
\begin{prop}
The system~\eqref{sysA14} has fifteen families of formal solutions
\(f = (f_0,\ldots,f_4)\) with simple pole and
the types of the formal solutions are determined by~\(\Res f\).
\end{prop}

The fifteen families of formal meromorphic solutions are divided into five
4-parameter families and ten 3-parameter families.  We further
classify the families from the viewpoint of the actions of B\"acklund
transformations.

Let \(f\)~be the formal solution with \(\Res f = (-1,-3,0,3,1)\), i.e.,
of case~(iii).  Substituting it into~\(g_i = s_2(f_i)\) (\(i=0,\ldots,4\)),
namely into
\[
  g_0 = f_0,\quad g_1 = f_1 - \frac{\alpha_2}{f_2},\quad g_2 = f_2,\quad
  g_3 = f_3 + \frac{\alpha_2}{f_2},\quad g_4 = f_4,
\]
we obtain the series expansion
\begin{align*}
  g_0 &= -\frac{1}{T} + O(1), &
  g_1 &= -(\alpha_1 + \alpha_2) T + O(T^2), & \\
  g_2 &= -\frac{1}{3}\alpha_2 T + O(T^3), &
  g_3 &= -(\alpha_3 + \alpha_2) T + O(T^2), &
  g_4  = \frac{1}{T} + O(1),
\end{align*}
which is the formal solution
with \(\Res g = (-1,0,0,0,1)\) of the system~\eqref{sysA14} with parameter
\(s_2(\alpha) = (s_2(\alpha_0),\ldots,s_2(\alpha_4))\).  Note
that this does not hold in the case of~\(\alpha_2 = 0\), since
\(g = f\) in that case.  Let us denote this fact as
\[
  \Res f = (-1,-3,0,3,1) \Rightarrow \Res s_2(f) = (-1,0,0,0,1)
  \quad\text{if } \alpha_2 \neq 0.
\]
Similarly, we obtain
\begin{align*}
  \Res f = (-1,0,0,0,1) & \Rightarrow \Res s_3(f) = (-1,0,1,0,0)
  \quad\text{if } \alpha_3 \neq 0, \\
  \Res f = (-1,0,0,0,1) & \Rightarrow \Res s_1(f) = (0,0,-1,0,1)
  \quad\text{if } \alpha_1 \neq 0, \\
  \Res f = (-1,0,1,0,0) & \Rightarrow \Res s_1(f) = (0,0,0,0,0)
  \quad\text{if } \alpha_1 \neq 0.
\end{align*}
Here \(\Res s_1(f) = (0,0,0,0,0)\) means that \(s_1(f)\)~is a formal
holomorphic solution.  Hence, generically, each formal
meromorphic solution can be transformed into a holomorphic
solution of the system~\eqref{sysA14} with different parameter by an
appropriate B\"acklund transformation.  Therefore it is convenient to
distinguish each family of formal meromorphic solutions
assigning to it a B\"acklund
transformation expressed as a product of~\(s_i\) (\(i=0,\ldots,4\)).
For instance, we assign B\"acklund transformations \(s_1\),~\(s_3 s_1\)
and~\(s_2 s_3 s_1\) to the families of formal meromorphic solutions with
\(\Res f = (-1,0,1,0,0)\),~\((-1,0,0,0,1)\) and~\((-1,-3,0,3,1)\),
respectively.  We also say simply that they are of type \((1)\),~\((13)\)
and~\((132)\), respectively.  Such classification of all the families of
formal meromorphic solutions of~\eqref{sysA14}
is given in Table~\ref{tab:res4}, where the family of holomorphic solutions
is denoted as type~\((\void)\).  The position~\(t_0\) of pole
is not counted into the arbitrary constants as before.

\begin{table}
\caption{Classification of the families of solutions
         of the system \protect\(\sysAone{4}\protect\)}
\label{tab:res4}
\[\begin{array}{c|c|c|c}
\hline
\text{type}&\Res f &\text{corresponding BT}
  &\text{\# of arbitrary constants} \\ \hline

(\void)& (0,0,0,0,0) &   \text{id} & 5 \\ \hline

(1)   & (-1,0,1,0,0) &         s_1 &   \\
(2)   & (0,-1,0,1,0) &         s_2 &   \\
(3)   & (0,0,-1,0,1) &         s_3 & 4 \\
(4)   & (1,0,0,-1,0) &         s_4 &   \\
(0)   & (0,1,0,0,-1) &         s_0 &   \\ \hline
(13)  & (-1,0,0,0,1) &     s_3 s_1 &   \\
(24)  & (1,-1,0,0,0) &     s_4 s_2 &   \\
(30)  & (0,1,-1,0,0) &     s_0 s_3 & 3 \\
(41)  & (0,0,1,-1,0) &     s_1 s_4 &   \\
(02)  & (0,0,0,1,-1) &     s_2 s_0 &   \\ \hline
(132) & (-1,-3,0,3,1)& s_2 s_3 s_1 &   \\
(243) & (1,-1,-3,0,3)& s_3 s_4 s_2 &   \\
(304) & (3,1,-1,-3,0)& s_4 s_0 s_3 & 3 \\
(410) & (0,3,1,-1,-3)& s_0 s_1 s_4 &   \\
(021) & (-3,0,3,1,-1)& s_1 s_2 s_0 &   \\
\hline
\end{array}\]
\end{table}

\section{Coordinates for formal meromorphic solutions
         of the system \protect\(\sysAone{4}\protect\)}
\label{sec:sys4-2}

Now we shall prove the convergence of the formal meromorphic solutions
obtained in the preceding section.  This section is devoted to the proof of
the following proposition.

\begin{prop}
Any formal meromorphic solution of the system \eqref{sysA14} converges.
\end{prop}

The idea of the proof is almost the same as the case of the system of
type~\(\typeAone{2}\) described in Subsection~\ref{sec:sys2-2}:
we will introduce appropriate change of variables
so that the system~\eqref{sysA14} is
transformed to a holomorphic system and the arbitrary constants contained
in the formal meromorphic solution of~\eqref{sysA14} are considered as a
initial condition of the holomorphic solution of the converted system.

While the system~\eqref{sysA14} has the fifteen families of formal
meromorphic solutions, we discuss the three typical families of type
\((1)\),~\((13)\) and~\((132)\) in sequence.  The
convergence of the solutions of the remaining types can be shown by the
rotations of indices.

\subsection{The family of type (1)}

Let \(f = (f_0,\ldots,f_4)\) be the formal meromorphic solution of type
\((1)\). Then the series expansion of~\(f\) is the form
of~\eqref{seriestype1}, and \(f_0\)~has simple pole at~\(t = t_0\) while
\(f_1\)~has a simple zero.  We note that \(1/f_0\)~has a simple zero
and \(f_0 f_1\)~has no pole:
\[
  f_0 f_1 = \alpha_1 - \bigpr{c_2^1 + \alpha_1 c_0^0}\, T + O(T^2).
\]
Since \(\alpha_1 - f_0 f_1\) has a simple zero, multiplying \(f_0\),
we have the expansion
\[
  f_0 (\alpha_1 - f_0 f_1) = \bigpr{c_2^1 + \alpha_1 c_0^0} + O(T),
\]
which has no terms of negative power of~\(T\) but has the constant term
containing the arbitrary constant~\(c_2^1\).
We also note that \(f_0 + f_2\) has no pole,
even though \(f_2\)~has a simple pole.
Therefore we introduce the following transformation:
\begin{equation}\label{patchtype1}
\begin{aligned}
        u_0 &= 1/f_0, \\
        u_1 &= f_0 (\alpha_1 - f_0 f_1), \\
  u_0 + u_2 &= f_0 + f_2, \\
        u_3 &= f_3, \\
  u_0 + \cdots + u_4 &= f_0 + \cdots + f_4.
\end{aligned}
\end{equation}
We can verify that this biholomorphic mapping
from~\(\BC^5 \setminus \{f_0 = 0\}\) to~\(\BC^5 \setminus \{u_0 = 0\}\)
transforms the system~\eqref{sysA14} to the system
\begin{equation}\begin{aligned}\label{sysA14_1}
u_0' &=
  2 u_0^3 u_1 - (\alpha_0 + 2 \alpha_1 - 1) u_0^2
  + (u_1 + u_2 - u_3 + u_4) u_0 - 1, \\
u_1' &=
  - 3 u_0^2 u_1^2 + (2 \alpha_0 + 4 \alpha_1 - 1) u_0 u_1
  - (u_1 + u_2 - u_3 + u_4) u_1 - (\alpha_0 + \alpha_1) \alpha_1,  \\
u_2' &=
  - 2 u_0^3 u_1 + (\alpha_0 + 2 \alpha_1 -1) u_0^2
  - (u_1 + u_2 - u_3 + u_4) u_0 \\ & \quad
  + (u_3 - u_4 - u_1) (u_0 + u_2)
  - 2 u_0 u_1 + \alpha_0 + 2\alpha_1 + \alpha_2 + 1,  \\
u_3' &=
  (u_4 - u_0 + u_1 - u_2) u_3 + \alpha_3, \\
u_4' &=
  3 u_0^2 u_1^2 - (2\alpha_0 + 4 \alpha_1 - 1) u_0 u_1
  + (u_1 + u_2 - u_3 + u_4) u_1 + (\alpha_0 + \alpha_1) \alpha_1 \\ & \quad
  + (u_4 + u_1) (u_0 + u_2 - u_3) + 2 u_0 u_1 + \alpha_4 - \alpha_1.
\end{aligned}\end{equation}
Here notice that
\( u_0'+\cdots+u_4' = \fivealphas{}{+}{+}{+}{+} \).

Since the right-hand sides of the above equations are
polynomials in~\(u = (u_0,\ldots,u_4)\), there exists
a unique holomorphic solution~\(u = u(t)\) with the initial condition
\(u(t_0) = (u_0(t_0),\ldots,u_4(t_0)) = (0, h_1, h_2, h_3, h_4)\),
where \(h_1, h_2, h_3, h_4\) are arbitrary complex numbers.
Let \(F\)~be a polynomial of~\(u\) such that \( F u_0 - 1 \)~is the
right-hand side of the first equation in~\eqref{sysA14_1}, namely,
\(u_0' = F u_0 - 1\).  Then it follows that
\(F\vert_{u_0 = 0} = u_1 + u_2 - u_3 + u_4\), and hence, setting
\(k = F(t_0) = F(u(t_0)) = h_1 + h_2 - h_3 + h_4,\) we have
\[
  u_0' (t_0) = F (t_0)\, u_0(t_0) - 1 = -1, \quad
  u_0''(t_0) = F'(t_0)\, u_0(t_0) + F(t_0)\, u_0'(t_0) = k.
\]
Therefore the Taylor expansion of~\(u_0=u_0(t)\) is
\begin{align*}
  u_0 = u_0(t) &= 0 + u_0'(t_0)\,T + \tfrac{1}{2}u_0''(t_0)\,T^2 + O(T^3) \\
               &= -T \,\bigl(1 + \tfrac{1}{2}k\,T + O(T^2)\bigr).
\end{align*}
We also have the Taylor expansion of \(u_1,\ldots,u_4\) as
\[
  u_1 = h_1 + O(T),\quad u_2 = h_2 + O(T),\quad
  u_3 = h_3 + O(T),\quad u_4 = h_4 + O(T),
\]
and, from~\eqref{patchtype1},
we see that the transform~\(f = (f_0,\ldots,f_4)\)
of~\(u = (u_0,\ldots,u_4)\) is of the form
\begin{align*}
f_0 &= 1 / u_0
     = -\tfrac{1}{T} + \tfrac{1}{2}k + O(T), \\
f_1 &= u_0 (\alpha_1 - u_0 u_1) \\
    &= -T \,\Bigl( 1 + \tfrac{1}{2}k\,T + O(T^2) \Bigr)
         \cdot \Bigl( \alpha_1 + h_1 T + O(T^2) \Bigr) \\
    &= -\alpha_1 T - \bigl(h_1 + \tfrac{1}{2}\alpha_1 k\bigr)\,T^2
       + O(T^3), \\
f_2 &= u_0 + u_2 - 1 / u_0
     = \tfrac{1}{T} + \bigl(h_2 - \tfrac{1}{2}k\bigr) + O(T), \\
f_3 &= u_3
     = h_3 + O(T), \\
f_4 &= (u_1 + u_4) - f_1
     = h_1 + h_4 + O(T).
\end{align*}
This expansion of~\(f\) coincides formally with the formal
meromorphic solution~\eqref{seriestype1} of~\eqref{sysA14} under the
change of the constants
\[
  c_0^0 = \tfrac{1}{2}k, \quad
  c_2^1 = - h_1 - \tfrac{1}{2}\alpha_1 k, \quad
  c_0^2 = h_2 - \tfrac{1}{2}k, \quad
  c_0^3 = h_3, \quad
  c_0^4 = h_1 + h_4 = k - h_2 + h_3.
\]
Hence the formal meromorphic solution~\eqref{seriestype1} must converge.

\subsection{The family of type (13)}

First, transform the series~\eqref{seriestype13} of type~(13) formally
by~\eqref{patchtype1}, then we have
\begin{align*}
        u_0 &= -T - c_0^0\,T^2 + O(T^3), \\
        u_1 &= -\bigpr{c_2^1 + \alpha_1 c_0^0} + O(T), \\
  u_0 + u_2 &= -\tfrac{1}{T} + c_0^0 + O(T), \\
        u_3 &= \phantom{-\tfrac{1}{T} + c_0^0}
               -\alpha_3 T + c_2^3\,T^2 + O(T^3), \\
  u_1 + u_4 &= \phantom{-}\tfrac{1}{T} + c_0^0 + O(T).
\end{align*}
The variable~\(u_1\) already contains the arbitrary constant \(c_2^1\) in its
constant term, so we need to pull out the~\(c_2^3\) appearing in~\(u_3\).
Noting that \(u_0 + u_2\)~has a simple pole while \(u_3\)~has a simple zero,
we construct Taylor expansions successively as follows:
\begin{gather*}
  (u_0+u_2) u_3 = \alpha_3 - \bigpr{c_2^3 + \alpha_3 c_0^0}\,T + O(T^2), \\
  (u_0+u_2)(\alpha_3 - (u_0+u_2)u_3) = -\bigpr{c_2^3 + \alpha_3 c_0^0} + O(T).
\end{gather*}
Hence we introduce the following transformation:
\begin{equation}\begin{aligned}\label{patchtype13}
        v_0 &= u_0, \\
        v_1 &= u_1, \\
  v_0 + v_2 &= 1 / (u_0 + u_2), \\
        v_3 &= (u_0 + u_2)(\alpha_3 - (u_0 + u_2) u_3), \\
  v_0 + \cdots + v_4 &= u_0 + \cdots + u_4.
\end{aligned}\end{equation}
Then we can verify that the transformation from~\(f = (f_0,\ldots,f_4)\)
to~\(v = (v_0,\ldots,v_4)\) is expressed
by
\begin{equation}\begin{aligned}\label{patchtype13f}
        v_0 &= 1 / f_0, \\
        v_1 &= f_0 f_1(\alpha_1 - f_0 f_1), \\
  v_0 + v_2 &= 1 / (f_0 + f_2), \\
        v_3 &= (f_0 + f_2)(\alpha_3 - (f_0 + f_2) f_3), \\
  v_0 + \cdots + v_4 &= f_0 + \cdots + f_4,
\end{aligned}\end{equation}
and the variable \(v\) satisfies the system
\begin{align*}
  v_0'&=
      2 v_0^3 v_1 - ( \alpha_0 + 2 \alpha_1) v_0^2
    + 2 v_0 (v_0 + v_2) ((v_0 + v_2) v_3 - \alpha_3) \\ & \quad
    + v_0 (v_0 + v_2) + v_0 (v_1 + v_3 + v_4) - 1, \\
  v_1' &=
    - 3 v_0^2 v_1^2 + 2 (\alpha_0 + 2 \alpha_1) v_0 v_1
    - 2 v_1 (v_0 + v_2) ((v_0 + v_2) v_3 - \alpha_3) \\ & \quad
    - v_1 (v_0 + v_2) - v_1 (v_1 + v_3 + v_4)
    - (\alpha_0 + \alpha_1) \alpha_1,  \\
  v_0' + v_2' &=
      2 (v_0 + v_2)^2 ( 2 v_0 v_1 - \alpha_0 - 2\alpha_1
    + 2 (v_0 + v_2) v_3 - \alpha_2 - 2\alpha_3) \\ & \quad
    + (v_0 + v_2)^2 + (v_1 + v_3 + v_4) (v_0 + v_2) - 1, \\
  v_3'&=
    - 3 (v_0 + v_2)^2 v_3^2 - 4 v_0 v_1 (v_0 + v_2) v_3
    + 2 (\alpha_0 + 2\alpha_1 + \alpha_2 + 2\alpha_3)(v_0 + v_2)v_3 \\ & \quad
    - (v_0 + v_2) v_3 - ( v_1 +  v_3 +  v_4 ) v_3
    + 2\alpha_3 v_0 v_1
    - (\alpha_0 + 2\alpha_1 + \alpha_2 + \alpha_3) \alpha_3, \\
  v_0' + \cdots + v_4' &= \fivealphas{}{+}{+}{+}{+},
\end{align*}
from which it follows that
\begin{align*}
v_2' &=
    2 (v_0 + v_2) ( (v_0 + v_2) v_3 - \alpha_3 ) v_2 \\ & \quad
  + (2 v_0 + v_2) ( 2 (v_0 v_1 - \alpha_1) - (\alpha_0 + \alpha_2) ) v_2
      + (v_0 + v_2) v_2 \\ & \quad
  + (v_1 + v_3 + v_4) v_2 - \alpha_0 v_0^2.
\end{align*}

Now let \(v = v(t)\)~be the unique holomorphic solution with
\[
  v(t_0) = (v_0(t_0), \ldots, v_4(t_0)) = (0, h_1, 0, h_3, h_4).
\]
Then we have
\[
  v_0'(t_0) = -1, \quad v_0''(t_0) = -k; \quad
  v_2'(t_0) = v_2''(t_0) = 0, \quad v_2^{(3)}(t_0) = -2\alpha_2, \quad
  v_2^{(4)} = -8\alpha_2 k,
\]
where \(k := h_1 + h_3 + h_4\).
This gives the Taylor expansion of~\(v\) as follows:
\begin{align*}
  v_0 &= -T \,\bigl(1 + \tfrac{1}{2}k\,T + O(T^2)\bigr), &
  v_1 &= h_1 + O(T),  & \\
  v_2 &= -\tfrac{1}{3}\alpha_2 T^3 (1 + k\,T) + O(T^5), &
  v_3 &= h_3 + O(T),  & v_4 (t) = h_4 + O(T).
\end{align*}
Noting that
\[
  v_0 + v_2 = -T \,\bigl(1 + \tfrac{1}{2}k\,T + O(T^2)\bigr),
\]
we have
\begin{align*}
f_0 &= 1 / v_0
     =-\tfrac{1}{T} + \tfrac{1}{2}k + O(T), \\
f_1 &= v_0 (\alpha_1 - v_0 v_1)
     =-\alpha_1 T - \bigl(h_1 + \tfrac{1}{2}\alpha_1 k\bigr)\,T^2 + O(T^3), \\
f_2 &= \frac{1}{v_0 + v_2}-\frac{1}{v_0}
     = \frac{-v_2}{v_0 (v_0 + v_2)}
     = \tfrac{1}{3}\alpha_2 T + 0 \cdot T^2 + O(T^3), \\
f_3 &= (v_0 + v_2) (\alpha_3 - (v_0 + v_2) v_3)
     =-\alpha_3 T - \bigl(h_3 + \tfrac{1}{2}\alpha_3 k\bigr)\,T^2 + O(T^3), \\
f_4 &= (v_0+\cdots+v_4) - (f_1+\cdots+f_4)
     = \tfrac{1}{T} + \tfrac{1}{2}k + O(T).
\end{align*}
The formal meromorphic solution~\eqref{seriestype13} formally
coincides with this meromorphic solutions if
\[
  c_0^0 = \tfrac{1}{2}k, \quad
  c_2^1 = - h_1 - \tfrac{1}{2}\alpha_1 k, \quad
  c_2^3 = - h_3 - \tfrac{1}{2}\alpha_3 k,
\]
and hence the formal solution \eqref{seriestype13} must be convergent.

\subsection{The family of type (132)}

We need transform the formal meromorphic solutions~\eqref{seriestype132} by
the change of variables~\eqref{patchtype13f} and observe the arbitrary
constants \(c_0^0\),~\(c_2^3\) and~\(c_4^2\),
which will appear in the expansion of~\(v = (v_0,\ldots,v_4)\).

First of all, we observe that
\begin{equation*}
v_0 = 1/f_0
    = -T - c_0^0\,T^2
      + \tfrac{1}{3}\Bigbk{\fivealphas{2}{+3}{+5}{+3}{+}+2\bigpr{c_0^0}^2}\,T^3
      + O(T^4),
\end{equation*}
and \(v_0\)~has a simple zero.  Secondly, from
\begin{align*}
v_3 &= (f_0 + f_2) (\alpha_3 - (f_0 + f_2) f_3) \\
    &= -\tfrac{3}{T^3} \Bigl( 1 - 2 c_0^0 T
       - \tfrac{1}{5}
         \Bigbk{\fivealphas{7}{+11}{+15}{+9}{+3} - 9 \bigpr{c_0^0}^2}\,T^2 \\
    &\qquad\quad
       - \Bigbk{\tfrac{1}{15} (\fivealphas{22}{+36}{+50}{+29}{+8})c_0^0
                - \tfrac{14}{5}\bigpr{c_0^0}^3 + 7 c_2^3}\,T^3 \Bigr) + O(T),
\end{align*}
it follows that \(1/v_3\)~has a zero of order three.  We will select
\(v_0\) and~\(1/v_3\) as coordinates with initial value~\(0\).
Since
\begin{align*}
v_1 + v_3
  &= f_0 (\alpha_1 - f_0 f_1)
       + (f_0 + f_2) (\alpha_3 - (f_0 + f_2) f_3) \\
  &= f_0 (\alpha_1 + \alpha_3 - f_0 (f_1 + f_3) - 2 f_2 f_3 )
       + f_2 (\alpha_3 - f_2 f_3) \\
  &= \bigl[\tfrac{1}{2}(\fivealphas{}{+}{+}{+}{+})c_0^0 - 2 c_2^3\bigr] + O(T),
\end{align*}
the value of the arbitrary constant~\(c_2^3\) can be recovered
from the initial value of~\(v_1 + v_3\).
For the remaining arbitrary constant~\(c_4^2\), we observe that
\begin{align*}
v_2
  &= \frac{1}{f_0 + f_2} - \frac{1}{f_0}
   = \frac{-f_2}{f_0 (f_0 + f_2)} \\
  &= \tfrac{1}{3}\alpha_2 T^3 + \tfrac{2}{3}\alpha_2 c_0^0\,T^4
     + \tfrac{1}{15}\alpha_2
       \Bigbk{\fivealphas{7}{+11}{+15}{+9}{+3}
              + 11\bigpr{c_0^0}^2}\,T^5 \\ & \quad
     + \Bigbk{\tfrac{1}{45}\alpha_2
              \bigpr{(\fivealphas{62}{+96}{+135}{+84}{+28}) c_0^0
                     + 26\bigpr{c_0^0}^3 - 30 c_2^3} - c_4^2}\,T^6 \\ & \quad
     + O(T^7),
\end{align*}
and then write it as
\begin{align*}
v_2 + c_4^2\,T^6
  &= \tfrac{1}{3}\alpha_2 T^3
     \Bigl(1 + 2 c_0^0 T + \tfrac{1}{5}
           \Bigbk{\fivealphas{7}{+11}{+15}{+9}{+3}
                  + 11\bigpr{c_0^0}^2}\,T^2 \\ & \qquad\qquad
           + \tfrac{1}{15}
             \Bigbk{(\fivealphas{62}{+96}{+135}{+84}{+28})c_0^0
                    + 26\bigpr{c_0^0}^3 - 30 c_2^3}\,T^3\Bigr) \\ & \qquad
     + O(T^7).
\end{align*}
Then, multiplying \(v_3\), we obtain
\begin{align*}
(v_2 + c_4^2\,T^6) v_3
  &= - \alpha_2 \Bigl( 1 + 0\cdot T + 0\cdot T^2
     + \tfrac{1}{3}\Bigbk{c_2^3 + (\alpha_2+\alpha_3) c_0^0}\,T^3 \Bigr)
     + O(T^4), \\
  &= - \alpha_2
     - \tfrac{1}{3}\alpha_2\Bigbk{c_2^3 + (\alpha_2+\alpha_3) c_0^0}\,T^3
     + O(T^4).
\end{align*}
Noting that
\(c_4^2\,T^6 \cdot v_3 = -3 c_4^2\,T^3 + O(T^4)\),
we have
\begin{align*}
  v_2 v_3 = -\alpha_2 -
    \Bigbk{&\tfrac{1}{3}\alpha_2
    \bigpr{c_2^3 + (\alpha_2+\alpha_3) c_0^0} - 3 c_4^2}\,T^3 + O(T^4),
\end{align*}
and hence \(-\alpha_2 - v_2 v_3\) has a zero of order three.

Thus we are now prepared to define the following transformation:
\begin{equation}\begin{aligned}\label{patchtype132}
        w_0 &= v_0, \\
  w_1 + w_3 &= v_1 + v_3, \\
        w_2 &= v_3 (-\alpha_2 - v_2 v_3), \\
        w_3 &= 1 / v_3,\\
  w_0 + \cdots + w_4 &= v_0 + \cdots + v_4.
\end{aligned}\end{equation}
We can verify that the system for
\(w\) is of the following holomorphic form:
\begin{align*}
w_0' &=
  2 w_0^3 (w_1 + w_3) - (\alpha_0 + 2 \alpha_1) w_0^2
  - 4 w_0^2 w_2 w_3 - 2 (2 \alpha_2 + \alpha_3) w_0^2 \\&\quad
  + 2 w_0 w_3 (w_2 w_3 + \alpha_2) (w_2 w_3 + \alpha_2 + \alpha_3) \\&\quad
  + w_0 (w_0 + w_2 + w_4) + w_0 (w_1 + w_3) - 1,
\\
w_1' + w_3' &=
  - 2 (w_1 + w_3) w_3 (w_2 w_3+\alpha_2)(w_2 w_3+\alpha_2+\alpha_3) \\&\quad
  - 3 w_0^2 (w_1 + w_3)^2 + 8 w_0 (w_1 + w_3) w_2 w_3 - w_2^2 w_3^2 \\&\quad
  + 2 (\alpha_0 + 2 \alpha_1 + 4 \alpha_2 + 2 \alpha_3) w_0 (w_1 + w_3)
  - 2 (\alpha_0 + 2 \alpha_1 + 2 \alpha_2 + \alpha_3) w_2 w_3 \\&\quad
  - 2 w_0 w_2 - (w_0 + w_2 + w_4) (w_1 + w_3) -  (w_1 + w_3)^2 \\&\quad
  - (\alpha_0 + \alpha_1 + 2 \alpha_2 + \alpha_3)
    (\alpha_1 + 2 \alpha_2 + \alpha_3)
  + \alpha_2 (\alpha_2 + \alpha_3),
\\
w_2' &=
  - w_2 w_3 (w_2 w_3 + \alpha_2) (w_2 w_3 + \alpha_2 + \alpha_3) \\&\quad
  + [ w_2 w_3 (w_2 w_3 + \alpha_2)
           + w_2 w_3 (w_2 w_3 + \alpha_2 + \alpha_3) \\&\qquad\quad
      + (w_2 w_3 + \alpha_2) (w_2 w_3 + \alpha_2 + \alpha_3)
      - 2 w_0 w_2 ] \\ &\qquad
    \times
    [ 2 w_0 (w_1 + w_3) - w_2 w_3
      - (\alpha_0 + 2 \alpha_1 + 2 \alpha_2 + \alpha_3) ] \\&\quad
   - (w_0 + w_2 + w_4) w_2 - (w_1 + w_3) w_2,
\\
w_3' &=
  w_3^2 (w_2 w_3 + \alpha_2) (w_2 w_3 + \alpha_2 + \alpha_3) \\&\quad
  - [ w_3^2 (2 w_2 w_3 + 2\alpha_2 + \alpha_3) - 2 w_0 w_3 ] \\&\qquad
    \times
    [ 2 w_0 (w_1 + w_3) - w_2 w_3
      - (\alpha_0 + 2 \alpha_1 + 2 \alpha_2 + \alpha_3) ] \\&\quad
  + (w_0 + w_2 + w_4) w_3+ (w_1 + w_3) w_3 - w_0^2,
\\
w_0' + \cdots + w_4' &= \fivealphas{}{+}{+}{+}{+}.
\end{align*}

Let \(w=w(t)\)~be a unique holomorphic solution of the above system with
\[
  w(t_0) = (w_0(t_0),\ldots,w_4(t_0)) = (0,h_1,h_2,0,h_4).
\]
Then, by a tedious calculation, we can obtain the expansion of~\(w\)
as
\begin{align*}
w_0 &= -T \,\bigl(1 + \tfrac{1}{2}k\,T + a_2\,T^2 + a_3\,T^3 + O(T^4)\bigr), &
w_1 &= h_1 + O(T), &
w_2 &= h_2 + O(T), & \\
w_3 &= -\tfrac{1}{3}T^3 \bigl(1 + k\,T + b_2\,T^2 + b_3\,T^3 + O(T^4)\bigr), &
w_4 &= h_4 + O(T), &
\end{align*}
where
\begin{align*}
 k  &= h_1 + h_2 + h_4, \\
a_2 &= \tfrac{1}{6}k^2 + \tfrac{1}{3}(\fivealphas{2}{+3}{+5}{+3}{+}), \\
a_3 &= \tfrac{1}{24}k^3 + \tfrac{1}{2}h_1
         + \tfrac{1}{24}(\fivealphas{13}{+21}{+37}{+21}{+5})k, \\
b_2 &= \tfrac{11}{20}k^2 + \tfrac{1}{5}(\fivealphas{7}{+11}{+15}{+9}{+3}), \\
b_3 &= \tfrac{13}{60}k^3 + \tfrac{7}{6}h_1
         + \tfrac{1}{120}(\fivealphas{213}{+349}{+485}{+281}{+77})k .
\end{align*}
From \eqref{patchtype1},~\eqref{patchtype13} and~\eqref{patchtype132},
it follows that
\begin{align*}
f_0 &= 1 / w_0
     = - \tfrac{1}{T} + \tfrac{1}{2}k + O(T),
\\
f_1 &= w_0(\alpha_1 - w_0 w_1) - w_0^2 w_3 + {w_0^2}/{w_3} \\
    &= - \tfrac{3}{T} + 0
       - \Bigl[\tfrac{1}{5}(\fivealphas{-}{+2}{+5}{+3}{+})
                + \tfrac{1}{10} k^2\Bigr]\,T \\&\quad
       - \Bigl[\tfrac{1}{2}h_1
                + \tfrac{1}{8}k(\fivealphas{}{+5}{+9}{+5}{+})\Bigr]\,T^2
       + O(T^3),
\\
f_2 &= \frac{1}{w_0-w_3(w_3w_2+\alpha_2)} - \frac{1}{w_0}
     = \frac{w_3(w_3w_2+\alpha_2)}{w_0(w_0-w_3(w_3w_2+\alpha_2))} \\
    &= - \tfrac{1}{3}\alpha_2 T + 0\cdot T^2
       - \tfrac{1}{3}\alpha_2
         \Bigl[\tfrac{1}{15}(\fivealphas{}{+3}{+0}{-3}{-4})
                - \tfrac{1}{30}k^2 \Bigr]\,T^3 \\&\quad
       + \tfrac{1}{9}
         \Bigl[ h_2 - \tfrac{1}{2}\alpha_2 h_1
           + \tfrac{1}{8}\alpha_2(\fivealphas{}{+}{+5}{+5}{+}) k\Bigr]\,T^4
       + O(T^5),
\\
f_3 &= - {w_0^2}/{w_3} + w_0 (2\alpha_2 + \alpha_3 + 2 w_2 w_3)
       - w_3 (\alpha_2 + w_2 w_3) (\alpha_2 + \alpha_3 + w_2 w_3) \\
    &= \tfrac{3}{T} + 0
       - \tfrac{1}{5}
         \Bigl[\fivealphas{}{+3}{+5}{+2}{-} - \tfrac{1}{2}k^2\Bigr]\,T \\&\quad
       + \Bigl[\tfrac{1}{8}(\fivealphas{}{+}{+}{+}{+})k
                - \tfrac{1}{2}h_1\Bigr]\,T^2
       + O(T^3),
\\
f_4 &= \tfrac{1}{T} + \tfrac{1}{2}k + O(T).
\end{align*}
Therefore, by the same argument as in the preceding subsections,
the formal meromorphic solutions~\eqref{seriestype132} converge.
The relations between the arbitrary constants in~\eqref{seriestype132}
and~\(h_1,h_2,h_4\) are given by
\begin{align*}
  c_0^0 &= \tfrac{1}{2}k, \quad
  c_2^3  =-\tfrac{1}{2}h_1 + \tfrac{1}{8}(\fivealphas{}{+}{+}{+}{+})k, \\
  c_4^2 &= \tfrac{1}{9}h_2 - \tfrac{1}{18}\alpha_2 h_1
           + \tfrac{1}{72}\alpha_2 (\fivealphas{}{+}{+5}{+5}{+})k.
\end{align*}

\section{Augmentation of the phase space
         of the system \protect\(\sysAone{4}\protect\)}
\label{sec:sys4-3}

We define a fiber space~\(\BE\) over the
parameter space
\[
  V = \condset{\alpha = (\alpha_0,\ldots,\alpha_4) \in \BC^5}
              {\alpha_0 + \cdots + \alpha_4 = 1},
\]
of which each fiber~\(\BE(\alpha)\) for~\(\alpha\in V\) is the augmented
phase space of the system~\eqref{sysA14}.  Then we observe how B\"acklund
transformations act on the space~\(\BE\).

First, we define the space~\(\BE\) and give its properties.  Let
\begin{align*}
  I = \bigl\{\void,
    0_+&,1_+,\ldots,4_+,02_+,13_+,\ldots,41_+,021_+,132_+,\ldots,410_+, \\
    0_-&,1_-,\ldots,4_-,02_-,13_-,\ldots,41_-,021_-,132_-,\ldots,410_- \bigr\},
\end{align*}
be a label set and let \(W_*\) (\(* \in I\)) be thirty-one copies of
\(V\times\BC^5\) with coordinate system
\((\alpha, x_*)=\bigpr{\alpha_0,\ldots,\alpha_4; x_*^0,\ldots,x_*^4} \in W_*\).
We will define identifying relations among \(W_*\).
In order to express identifying relations in simple form,
we use some auxiliary mappings.

Let
\(\Psi_{i_+},\Psi_{i_-}, \overline\Psi_{i_+}, \overline\Psi_{i_-}\)
(\(i = 0,\ldots,4\)) be birational mappings given by
\(\Psi_{i_\pm} (\alpha,\xi) = (\alpha,\eta)\) and
\(\overline\Psi_{i_\pm} (\alpha, \overline\xi) = (\alpha, \overline\eta)\)
where
\begin{gather*}
\begin{split}
  \eta^{i\pm 1} = 1/\xi^{i\pm 1}, \quad
& \eta^i \eta^{i\pm 1} = \mp \alpha_i - \xi^i \xi^{i\pm 1}, \\
& \eta^{i\pm 1} + \eta^{i\mp 1} = \xi^{i\pm 1} + \xi^{i\mp 1}, \quad
  \eta^{i\mp 2} = \xi^{i\mp 2},\ \sum_j\eta^j = \sum_j\xi^j;
\end{split} \\
\begin{split}
  \overline\eta^{i\pm 3} = \overline\xi^{i\pm 3}, \quad
& \overline\eta^{i\pm 2} = \overline\xi^{i\pm 2}, \quad
  \overline\eta^{i\pm 3} + \overline\eta^{i\pm 1}
    = 1/(\overline\xi^{i\pm 3} + \overline\xi^{i\pm 1}), \\
& \overline\eta^i (\overline\eta^{i\pm 3} + \overline\eta^{i\pm 1})
    = \mp\alpha_i -
      \overline\xi^i (\overline\xi^{i\pm 3} + \overline\xi^{i\pm 1}), \quad
  \sum_j\overline\eta^j = \sum_j\overline\xi^j.
\end{split}
\end{gather*}
Then we define the space~\(\BE\) by gluing~\(W_*\)
via the following identifying equations
among coordinates~\((\alpha, x_*)\) and~\((\alpha, x_\void)\)
\begin{gather*}
  \Psi_{i_+}(\alpha, x_\void) = (\alpha, x_{i_+}), \quad
  \overline\Psi_{i-1_+}(\alpha, x_{i+1_+}) = (\alpha, x_{i-1,i+1_+}), \quad
  \Psi_{i_+}(\alpha, x_{i-1,i+1_-}) = (\alpha, x_{i-1,i+1,i_+}), \\
  \Psi_{i_-}(\alpha, x_\void) = (\alpha, x_{i_-}), \quad
  \overline\Psi_{i+1_-}(\alpha, x_{i-1_-}) = (\alpha, x_{i-1,i+1_-}), \quad
  \Psi_{i_-}(\alpha, x_{i-1,i+1_+}) = (\alpha, x_{i-1,i+1,i_-}).
\end{gather*}
That is,
\[
  \BE = \biggl(\bigsqcup_{*\in I}W_*\biggr)\,\bigg/\sim,
\]
where \(\sim\)~is the equivalence relation generated by the above
equations.  Here, we notice that the coordinate systems in
\eqref{patchtype1},~\eqref{patchtype13} and~\eqref{patchtype132} are
written as
\[
  f = x_\void, \quad u = x_{1_-}, \quad v = x_{13_-}, \quad w = x_{132_+}
\]
in our new notation.

In the same way as in Subsection~\ref{sec:sys2-3},
we define the subset~\(U_*\) of~\(\BE\)
for each~\(* \in I\) and coordinate mappings
\[
  \phi_* \colon p \in U_* \mapsto
  \bigpr{\alpha_0(p), \ldots, \alpha_4(p); x_*^0(p), \ldots, x_*^4(p)}
    \in W_* = V\times\BC^5 \quad (*\in I).
\]
Then the space~\(\BE\) is described by the atlas~\(\{(\phi_*, U_*)\}\).
We also have the natural projection~\(\pi_V \colon \BE \to V\).
The fiber~\(\BE(\alpha) := \pi_V^{-1}(\alpha)\) is
a five dimensional complex manifold
for each~\(\alpha \in V\).

Now we see that the space~\(\BE\) is decomposed into the disjoint union of
subsets, each of which corresponds to a family of meromorphic solutions.
Let
\begin{align*}
  \BE_\void &= \Uvoid, \\
  \BE_1     &= \Usubset{1_-}   = \invzero{1_-}{0}{} \\
            &= \Usubset{1_+}   = \invzero{1_+}{2}{},\\
  \BE_{13}  &= \Usubset{13_-}  = \invzero{13_- }{0}{2} \\
            &= \Usubset{13_+}  = \invzero{13_+ }{4}{2},\\
  \BE_{132} &= \Usubset{132_-} = \invzero{132_-}{0}{3} \\
            &= \Usubset{132_+} = \invzero{132_+}{4}{1}.
\end{align*}
Then
\[
  \BE_\void \cong V \times \BC^5,
\ \BE_1     \cong V \times \BC^4,
\ \BE_{13}  \cong V \times \BC^3,
\ \BE_{132} \cong V \times \BC^3.
\]
Similarly, we define the following subspaces by cyclic rotation
\[
  \BE_2,     \BE_3,     \BE_4,     \BE_0;
  \BE_{24},  \BE_{30},  \BE_{41},  \BE_{02};
  \BE_{243}, \BE_{304}, \BE_{410}, \BE_{021}.
\]
Then the space~\(\BE\) is decomposed by these subsets as
\begin{align*}
\BE=\BE_\void&\sqcup \BE_0     \sqcup \BE_1     \sqcup\cdots\sqcup \BE_4    \\
             &\sqcup \BE_{02}  \sqcup \BE_{13}  \sqcup\cdots\sqcup \BE_{41} \\
             &\sqcup \BE_{021} \sqcup \BE_{132} \sqcup\cdots\sqcup \BE_{410}.
\end{align*}

The study of Section~\ref{sec:sys4-2} shows:
\begin{thm}
The vector field
\[
  X_p = \sum_{i=0}^4
        \Bigl(x_\void^i(p) \bigl(x_\void^{i+1}(p) - x_\void^{i+2}(p)
                               + x_\void^{i+3}(p) - x_\void^{i+4}(p)\bigr)
              + \alpha_i(p)\Bigr)
        \left(\frac\partial{\partial x_\void^i}\right)_p
\]
defined on~\(\BE(\alpha)\cap\Uvoid\) extends to the entire
fiber~\(\BE(\alpha)\) as a holomorphic vector field.
\end{thm}

Secondly, we investigate how B\"acklund transformations act
on the space~\(\BE\).  For any~\(w \in \WtypeAone{4}\),
we define a birational mapping~\(\sigma_w\) from~\(\BE\) to itself by
\[
  (\alpha_j\circ\sigma_w)(p) = (w(\alpha_j))(p), \quad
  \bigpr{x_*^j\circ\sigma_w}(p) = \bigpr{w\bigpr{x_*^j}}(p), \quad
    p\in\BE,\, j=0,\ldots,4
\]
for any~\(*\in I\).  Since \(\sigma_{w'}\circ\sigma_w = \sigma_{ww'}\)
for~\(w, w' \in\WtypeAone{4}\), and since \(\WtypeAone{4}\)~is generated
by~\(s_0,\ldots,s_4\) and~\(\pi\), we need only to study
\(\sigma_i := \sigma_{s_i}\) (\(i = 0,\ldots,4\)) and~\(\sigma_\pi\).

We first note that:
\begin{thm}
For any \(w \in \WtypeAone{4}\), \(\sigma_w\) is a biholomorphic mapping
from \(\BE(\alpha)\) to \(\BE(w(\alpha))\) for any \(\alpha\in V\), and the
vector field \(X\) is invariant by \(\sigma_w\).
\end{thm}
\begin{proof}
It is sufficient to show the proposition for \(w = s_0,\ldots,s_4,\pi\)
and these special cases can be shown by explicit calculation.  By use of
the coordinates \(u\),~\(v\) and~\(w\) defined in the preceding section,
we can verify that
\begin{align*}
  s_i(u_j) &=
    \begin{cases}
      u_j + \dfrac{\alpha_i}{u_i} u_{ij}, & i = 1, 3, \\
      u_j,                                & i = 0.
    \end{cases} & (u_j = x_{1_-}^j) & \\
  s_i(v_j) &=
    \begin{cases}
      v_j + \dfrac{\alpha_i}{v_i} u_{ij}, & i = 1, 2, 3, \\
      v_j,                                & i = 0.
    \end{cases} & (v_j = x_{13_-}^j) & \\
  s_i(w_j) &=
    \begin{cases}
      w_j + \dfrac{\alpha_i}{w_i} u_{ij}, & i = 2, \\
      w_j,                                & i = 0, 3.
    \end{cases} & (w_j = x_{132_+}^j) &
\end{align*}
The remaining coordinates can be computed similarly.
\end{proof}

\begin{rem}
For any~\(\alpha\in V\) and~\(w \in \WtypeAone{4}\), the augmented phase
spaces \(\BE(\alpha)\)~and~\(\BE(w(\alpha))\) are isomorphic.
\end{rem}

Now we study how the mapping~\(\sigma_i\) (\(i = 0,\ldots,4\)) acts on each
component of the decomposition of~\(\BE\).  We observe here only the case
of~\(i = 2\), since the other cases are obtained by cyclic rotations.  Let
\begin{align*}
  D_2 = \overline{\bigl\{x_\void^2 = 0\bigr\}}
      = \bigl\{x_\void^2 = 0\bigr\} \cup
        \bigl\{x_{0_-}^2 = 0\bigr\} \cup
        \bigl\{x_{4_+}^2 = 0\bigr\} \cup
        \bigl\{x_{13_-}^2 = 0\bigr\}
\end{align*}
and
\[
  \BE_*' = \BE_* \setminus \pi_V^{-1}(\{\alpha_2 = 0\}).
\]
Here \(\invzero{*}{2}{}\) is simply denoted
by~\(\bigl\{x_*^2=0\bigr\}\) and the closure is that in the space~\(\BE\).
Then the points of~\(\BE \setminus \pi_V^{-1}(\{\alpha_2 = 0\})\) are
mapped as follows:
\begin{gather*}
  \begin{aligned}
  \sigma_2(\BE_\void \setminus D_2) &= \BE_\void \setminus D_2, &
  \sigma_2(\BE_\void' \cap D_2) &= \BE_2', &
  \sigma_2(\BE_2') &= \BE_\void' \cap D_2, \\
  \sigma_2(\BE_0 \setminus D_2) &= \BE_0 \setminus D_2, &
  \sigma_2(\BE_0' \cap D_2) &= \BE_{02}',  &
  \sigma_2(\BE_{02}') &= \BE_0' \cap D_2, \\
  \sigma_2(\BE_4 \setminus D_2) &= \BE_4 \setminus D_2, &
  \sigma_2(\BE_4' \cap D_2) &= \BE_{24}', &
  \sigma_2(\BE_{24}') &= \BE_4' \cap D_2,
  \end{aligned} \\
  \sigma_2(\BE_{13}' \cap D_2) = \BE_{132}', \ \qquad
  \sigma_2(\BE_{132}') = \BE_{13}' \cap D_2, \\
  \sigma_2(\BE_1) = \BE_1, \quad
  \sigma_2(\BE_{41}) = \BE_{41}, \quad
  \sigma_2(\BE_{410}) = \BE_{410}, \quad
  \sigma_2(\BE_{021}) = \BE_{021}, \\
  \sigma_2(\BE_3) = \BE_3, \quad
  \sigma_2(\BE_{30}) = \BE_{30}, \quad
  \sigma_2(\BE_{304}) = \BE_{304}, \quad
  \sigma_2(\BE_{243}) = \BE_{243}.
\end{gather*}
These relationships explain the correspondence between the types of
meromorphic solutions and the B\"acklund transformations.

To end this section, we remark that the space~\(\BE\) contains the
augmented phase space for the system~\(\sysAone{2}\) as a submanifold with
an appropriate restriction.  In fact, setting~\(\alpha_3 = \alpha_4 = 0\),
the system~\eqref{sysA14} can be restricted to~\(f_3 = f_4 = 0\).
Therefore \(\overline{(\phi_\void)^{-1}
\bigl(\bigl\{\alpha_3 = \alpha_4 = x_\void^3 = x_\void^4 = 0\bigr\}\bigr)}
\subset\BE\) is isomorphic to the fiber space for the system~\(\sysAone{2}\).

\section{The \protect\(\typeAone{4}\protect\) Hamiltonian system}
\label{sec:sys4-hsys}

The differential system~\eqref{sysA14} is equivalent to the Hamiltonian
system, as is shown in~\cite{NY}.  This is shown by introducing a new
coordinate system~\((p_1, q_1, p_2, q_2; t)\) in the original
phase space~\(\condset{f\in\BC^5}{}\) fixing~\(f_0 + \cdots + f_4 = t\).
Then the system~\eqref{sysA14} is written as a Hamiltonian system
\begin{gather*}
  \frac{dq_i}{dt} =  \frac{\partial H}{\partial p_i},\quad
  \frac{dp_i}{dt} = -\frac{\partial H}{\partial q_i} \quad (i=1,2)
\end{gather*}
with
a polynomial Hamiltonian function \(H\) of \(p_1, q_1, p_2, q_2\) and~\(t\).
We can also choose a coordinate system of each~\(U_*\) for~\(*\in I\)
so that the transformation from the coordinate system of~\(\Uvoid\)
to that of~\(U_*\) is symplectic. This means that the Hamiltonian system
in~\(\Uvoid\) extends to the whole space~\(\BE\).

We list below the Hamiltonian functions~\(H_*\)
for \(* = \void,\ 1_+,\ 1_-,\ 13_+,\ 13_-,\ 132_+,\ 132_-\), where the
canonical coordinates in~\(U_*\) is denoted by the same notation
\((p_1, q_1, p_2, q_2)\) for simplicity.  These canonical coordinates are
written by the original coordinates~\(f\) of the system~\eqref{sysA14}
as well as by our new coordinates~\(x_*\).
Note that we fix \(x_*^0 + \cdots + x_*^4 = t\) in every case.

\begin{description}

\item [ The Hamiltonian \(H_\void\) in \(\Uvoid\) ]
\begin{align*}
H_\void &=
  (t - q_1 - p_1) q_1 p_1 + (t - q_2 - p_2) q_2 p_2 - 2 q_1 p_1 q_2 \\&\quad
  - \alpha_1 q_1 + \alpha_2 p_1 - (\alpha_1 + \alpha_3) q_2 + \alpha_4 p_2,
\end{align*}
where
\begin{gather*}
   p_1 = f_1,\ q_1 = f_2,\ p_2 = f_1 + f_3,\ q_2 = f_4;
\\
  (p_1, q_1, p_2, q_2) =
    (x_\void^1, x_\void^2, x_\void^1 + x_\void^3, x_\void^4).
\end{gather*}

\item [ The Hamiltonian \(H_{1_+}\) in \(U_{1_+}\) ]
\begin{align*}
H_{1_+} &=
  (t - q_1 - p_1) q_1 p_1 - t (q_2 p_2 + \alpha_1)
  - q_2 (q_2 p_2 + \alpha_1) (q_2 p_2 + \alpha_1 + \alpha_2) \\&\quad
  + 2 p_1 q_2 p_2 + (\alpha_0 + 2\alpha_1 + \alpha_2) p_1
  - \alpha_4 q_1 + p_2,
\end{align*}
where
\begin{gather*}
  p_1 = f_4,\ q_1 = f_0 + f_2,\ p_2 = f_2 (-\alpha_1 - f_1 f_2),\ q_2 = 1/f_2;
\\
  (p_1, q_1, p_2, q_2) =
    (x_{1_+}^4, x_{1_+}^0\!\!+ x_{1_+}^2, x_{1_+}^1, x_{1_+}^2).
\end{gather*}

\item [ The Hamiltonian \(H_{1_-}\) in \(U_{1_-}\) ]
\begin{align*}
H_{1_-} &=
  (t - q_2 - p_2) q_2 p_2 - t (q_1 p_1 - \alpha_1)
  - p_1 (\alpha_1 - q_1 p_1) (\alpha_0 + \alpha_1 - q_1 p_1) \\&\quad
  + 2 q_1 p_1 q_2 + q_1 - (\alpha_0 + 2\alpha_1 + \alpha_2) q_2 + \alpha_3 p_2,
\end{align*}
where
\begin{gather*}
   p_1 = 1/f_0,\ q_1 = f_0 (\alpha_1 - f_0 f_1),\ p_2 = f_0 + f_2,\ q_2 = f_3;
\\
  (p_1, q_1, p_2, q_2) =
    (x_{1_-}^0, x_{1_-}^1, x_{1_-}^0\!\!+ x_{1_-}^2, x_{1_-}^3).
\end{gather*}

\item [ The Hamiltonian \(H_{13_+}\) in \(U_{13_+}\) ]
\begin{align*}
H_{13_+} =&
  - q_1 (q_1 p_1 + \alpha_1) (q_1 p_1 + \alpha_1 + \alpha_2)
  - q_2 (q_2 p_2 + \alpha_3) (q_2 p_2 + \alpha_3 + \alpha_4) \\ &
  - q_1 (q_1 p_1 + \alpha_1) (2 q_2 p_2 + 2\alpha_3 + \alpha_4)
  - t (q_1 p_1 + \alpha_1 + q_2 p_2 + \alpha_3)
  + p_1 + p_2,
\end{align*}
where
\begin{gather*}
   p_1 = (f_2 + f_4) (-\alpha_1 - (f_2 + f_4) f_1),\ q_1 = 1/(f_2 + f_4),
 \ p_2 = f_4 (-\alpha_3 - f_3 f_4),\ q_2 = 1/f_4;
\\
  (p_1, q_1, p_2, q_2) =
    (x_{13_+}^1, x_{13_+}^2\!\!+ x_{13_+}^4, x_{13_+}^3, x_{13_+}^4).
\end{gather*}

\item [ The Hamiltonian \(H_{13_-}\) in \(U_{13_-}\) ]
\begin{align*}
H_{13_-} = &
  - p_1 (\alpha_1 - q_1 p_1) (\alpha_0 + \alpha_1 - q_1 p_1)
  - p_2 (\alpha_3 - q_2 p_2) (\alpha_2 + \alpha_3 - q_2 p_2) \\ &
  - p_2 (\alpha_3 - q_2 p_2) (\alpha_0 + 2\alpha_1 - 2 q_1 p_1)
  + t (\alpha_1 - q_1 p_1 + \alpha_3 - q_2 p_2 )
  + q_1 + q_2,
\end{align*}
where
\begin{gather*}
   p_1 = 1/f_0,\ q_1 = f_0 (\alpha_1 - f_0 f_1),
 \ p_2 = 1/(f_0 + f_2),\ q_2 = (f_0 + f_2) (\alpha_3 - (f_0 + f_2) f_3);
\\
  (p_1, q_1, p_2, q_2) =
    (x_{13_-}^0, x_{13_-}^1, x_{13_-}^0\!\!+ x_{13_-}^2, x_{13_-}^3).
\end{gather*}

\item [ The Hamiltonian \(H_{132_+}\) in \(U_{132_+}\) ]
\begin{align*}
H_{132_+} &=
  q_2 (q_2 p_2 + \alpha_2) (q_2 p_2 + \alpha_2 + \alpha_3)
  (\alpha_0 + 2\alpha_1 - 2 q_1 p_1 + q_2 p_2 + 2\alpha_2 + \alpha_3) \\&\quad
  - p_1 (\alpha_1 - q_1 p_1 + q_2 p_2 + 2\alpha_2 + \alpha_3)
    (\alpha_0 + \alpha_1 - q_1 p_1 + q_2 p_2 + 2\alpha_2 + \alpha_3) \\&\quad
  - p_1 (q_2 p_2 (\alpha_0 + 2\alpha_1 - 2 q_1 p_1)
         - \alpha_2 (\alpha_2 + \alpha_3) + p_1 p_2) \\ &\quad
  + t (\alpha_1 - q_1 p_1 + q_2 p_2  + \alpha_2 + \alpha_3)
  + q_1,
\end{align*}
where
\begin{gather*}
   p_1 = 1/f_0,\ q_1 - 1/q_2 = f_0 (\alpha_1 - f_0 f_1), \\
   q_2 (-\alpha_2 - q_2 p_2) + p_1 = 1/(f_0 + f_2),
 \ 1/q_2 = (f_0 + f_2)(\alpha_3 - (f_0 + f_2) f_3);
\\
  (p_1, q_1, p_2, q_2) =
    (x_{132_+}^0, x_{132_+}^1\!\!+ x_{132_+}^3, x_{132_+}^2, x_{132_+}^3).
\end{gather*}

\item [ The Hamiltonian \(H_{132_-}\) in \(U_{132_-}\) ]
\begin{align*}
H_{132_-} = &
  - p_1 (\alpha_2 - q_1 p_1) (\alpha_1 + \alpha_2 - q_1 p_1)
    (\alpha_1 + 2\alpha_2 - q_1 p_1 + 2 q_2 p_2 + 2\alpha_3 + \alpha_4) \\ &
  - q_2 (\alpha_1 + 2\alpha_2 - q_1 p_1 + q_2 p_2 + \alpha_3)
        (\alpha_1 + 2\alpha_2 - q_1 p_1 + q_2 p_2 + \alpha_3 + \alpha_4) \\ &
  + q_2 (q_1 p_1 (2 q_2 p_2 + 2\alpha_3 + \alpha_4)
         + (\alpha_1 + \alpha_2) \alpha_2 - q_1 q_2) \\ &
  - t (\alpha_1 + \alpha_2 - q_1 p_1 + q_2 p_2  + \alpha_3)
  + p_2,
\end{align*}
where
\begin{gather*}
   1/p_1 = (f_2 + f_4)(-\alpha_1 - (f_2 + f_4) f_1),
 \ p_1 (\alpha_2 - p_1 q_1) + q_2 = 1/(f_2 + f_4), \\
   p_2 - 1/p_1 = f_4 (-\alpha_3 - f_3 f_4),\ q_2 = 1/f_4;
\\
  (p_1, q_1, p_2, q_2) =
    (x_{132_-}^1, x_{132_-}^2, x_{132_-}^1\!\!+x_{132_-}^3, x_{132_-}^4).
\end{gather*}

\end{description}

\bigskip

\noindent
\textit{Acknowledgment.}\quad
The author express his sincere gratitude to
Professor Kyoichi Takano
who read this paper carefully and gave helpful advice and encouragement.

\bigskip

\begin{flushright}
\begin{tabular}[h]{c}
Nobuhiko Tahara \\
Graduate School of Science and Technology \\
Kobe University \\
Rokko, Kobe \\
657-8501 Japan \\
(E-mail: \texttt{tahara@math.kobe-u.ac.jp})
\end{tabular}
\end{flushright}

\end{document}